\documentclass[12pt]{article}
\usepackage{amsmath}
\usepackage{amssymb}
\usepackage{bm}
\usepackage{graphicx,psfrag,epsf}
\usepackage{enumerate}
\usepackage{url}
\usepackage{amssymb}
\usepackage{amsmath}
\usepackage{natbib}
\usepackage{amsthm}
\usepackage{booktabs}
\usepackage{arydshln}
\usepackage{mathtools}
\mathtoolsset{showonlyrefs,showmanualtags}
\usepackage[inline]{trackchanges}
\newcommand{\Var}{\mathrm{Var}} 
\newcommand{\Cov}{\mathrm{Cov}}
\newcommand{\Corr}{\mathrm{Corr}}

\begin{document}

%\bibliographystyle{natbib}

%Commands:
\def\spacingset#1{\renewcommand{\baselinestretch}%
{#1}\small\normalsize} \spacingset{1}

\newcommand{\indep}{\perp \!\!\! \perp}
\newcommand\independent{\protect\mathpalette{\protect\independenT}{\perp}}
\newcommand{\argmax}[1]{\underset{#1}{\operatorname{argmax}}\;} 
\newcommand{\argmin}[1]{\underset{#1}{\operatorname{argmin}}\;} 

\newtheorem{manualtheoreminner}{Theorem}
\newenvironment{manualtheorem}[1]{%
  \renewcommand\themanualtheoreminner{#1}%
  \manualtheoreminner
}{\endmanualtheoreminner}

\newtheorem{manuallemmainner}{Lemma}
\newenvironment{manuallemma}[1]{%
  \renewcommand\themanuallemmainner{#1}%
  \manuallemmainner
}{\endmanuallemmainner}

\newtheorem{definition}{Definition}
\newtheorem{result}{Result}
\newtheorem{example}{Example}
\newtheorem{theorem}{Theorem}
\newtheorem{corollary}{Corollary}
\newtheorem{lemma}{Lemma}
\newtheorem{proposition}{Proposition}
\newtheorem{assumption}{Assumption}

%Instructions:
%https://amstat.tfjournals.com/asa-style-guide/

%%%%%%%%%%%%%%%%%%%%%%%%%%%%%%%%%%%%%%%%%%%%%%%%%%%%%%%%%%%%%%%%%%%%%%%%%%%%%%

\title{\bf Adaptive Sequential Design for a Single Time-Series}
\author{
Ivana Malenica \\
Division of Biostatistics, University of California, Berkeley
\hspace{.2cm}\\
Aurelien Bibaut \\
Division of Biostatistics, University of California, Berkeley
\hspace{.2cm}\\
Mark J. van der Laan\\
Division of Biostatistics, University of California, Berkeley}

\date{}
\maketitle

\bigskip
\begin{abstract}
The current work is motivated by the need for robust statistical methods for precision medicine; we pioneer the concept of a sequential, adaptive design for a single individual. As such, we address the need for statistical methods that provide actionable inference for a single unit at any point in time. Consider the case that one observes a single time-series, where at each time $t$, one observes a data record $O(t)$ involving treatment nodes $A(t)$, an outcome node $Y(t)$, and time-varying covariates $W(t)$. We aim to learn an optimal, unknown choice of the controlled components of the design in order to optimize the expected outcome; with that, we adapt the randomization mechanism for future time-point experiments based on the data collected on the individual over time. Our results demonstrate that one can learn the optimal rule based on a single sample, and thereby adjust the design at any point $t$ with valid inference for the mean target parameter. This work provides several contributions to the field of statistical precision medicine. First, we define a general class of averages of conditional causal parameters defined by the current context (``context-specific'') for the single unit time-series data. 
%This estimation problem is addressed in a statistical model that is nonparametric under few assumptions on the data generating distribution. 
We define a nonparametric model for the probability distribution of the time-series under few assumptions, 
%that refrains from making unrealistic assumptions, 
and aim to fully utilize the sequential randomization in the estimation procedure via the double robust structure of the efficient influence curve of the proposed target parameter. We present multiple exploration-exploitation strategies for assigning treatment, and methods for estimating the optimal rule. Lastly, we present the study of the data-adaptive inference on the mean under the optimal treatment rule, where the target parameter adapts over time in response to the observed context of the individual. Our target parameter is pathwise differentiable with an efficient influence function that is doubly robust - which makes it easier to estimate than previously proposed variations. We characterize the limit distribution of our estimator under a Donsker condition expressed in terms of a notion of bracketing entropy adapted to martingale settings.
\end{abstract}

\noindent%
{\it Keywords:} Sequential decision making, Time-Series, Optimal Individualized Treatment, Targeted Maximum Likelihood Estimation (TMLE), Causal Inference.
\vfill

\newpage
\spacingset{1.45}
\section{Introduction}
\label{sec:intro}

There is growing scientific enthusiasm for the use, and development, of mobile health designs (mHealth) - broadly referring to the practice of health care mediated via mobile and wearable technologies \citep{steinhubl2013,malvey2014,istepanian2017,istepanian2018}. Numerous smartphones and Internet coupled devices, connected to a plethora of mobile health applications, support continuous assembly of data-driven healthcare intervention and insight opportunities. Interest in mobile interventions spans myriad of applications, including behavioral maintenance or change \citep{free2013,muhammad2017}, disease management \citep{heron2010, ertin2011AutoSenseUW, muessig2013, steinhubl2013}, teaching and social support \citep{kumar2013} and addiction management \citep{dulin2014,zhang2016tracker}. In particular, \citet{istepanian2018} refer to mHealth as one of the most transformative drivers for healthcare delivery in modern times. 

Recently, a new type of an experimental design termed micro-randomized trial (MRT) was developed in order to support just-in-time adaptive exposures - with an aim to deliver the intervention at the optimal time and location \citep{dempsey2015, klasnja2015}. To this date, multiple trials have been completed using MRT design, including encouraging regular physical activity \citep{klasnja2019} and engaging participation in substance use data gathering process in high-risk populations \citep{rabbi2018}. For both observational mHealth and MRT, the time-series nature of the collected data provides an unique opportunity to collect individual characteristics and context of each subject, while studying the effect of treatment on the outcome at specified future time-point.

The generalized estimating equation (GEE) and random effects models are the most commonly employed approaches for the analysis of mobile health data \citep{walls2006, stone2007, bolger2013intensive, hamaker2018}. As pointed out in \citet{boruvka2018}, these methods often do not yield consistent estimates of the causal effect of interest if time-varying treatment is present. As an alternative, \citet{boruvka2018} propose a centered and weighted least square estimation method for GEE that provides unbiased estimation, assuming linear model for the treatment effects. They tackle proximal and distal effects, with a focus on continuous outcome. On the other hand, \citet{luckett2019} propose a new reinforcement learning method applicable to perennial, frequently collected longitudinal data. While the literature on dynamic treatment regimes is vast and well-studied \citep{Murphy2003,Robins2004,Chakraborty2013,luedtke2016,luedtke2016resource,luedtke2016mean}, the unique challenges posed by mHealth obstruct their direct employment; for instance, mHealth objective typically has an infinite horizon. \citet{luckett2019} model the data-generating distribution as a Markov decision process, and estimate the optimal policy among a class of pre-specified policies in both offline and online setting. 

While mHealth, MRT designs and the corresponding methods for their analysis aim to deliver treatment tailored to each patient, they are still not optimized with complete ``N-of-1'' applications in mind. The usual population based target estimands fail to ensnare the full, personalized nature of the time-series trajectory, often imposing strong assumptions on the dynamics model for the estimation purposes. To the best of our knowledge, \citet{robins1999} provide the first step towards describing a causal framework for a single subject with time-varying exposure and binary outcome in a time-series setting. Focusing on full potential paths, \citet{bojinov2019} provide a causal framework for time-series experiments with randomization-based inference. Other methodologies focused on single unit applications rely on strong modeling assumptions, primarily linear predictive models and stationarity; see \citet{bojinov2019} for an excellent review of the few works on the topic.  %Overall, the focus in all is on full potential paths, not exploring the current context of the subject (effect moderation), and often depending on final time-point definition or complete randomization. 
Alternatively, \citet{vanderLaan2018onlinets} propose causal effects defined as marginal distributions of the outcome at a particular time point under a certain intervention on one or more of the treatment nodes. The efficient influence function of these estimators, however, relies on the whole mechanism in a non-double robust manner. Therefore, even when the assignment function is known, the inference still relies on consistent (at rate) estimation of the conditional distributions of the covariate and outcome nodes.

The current work presented is motivated by the need for robust statistical methods for precision medicine, pioneering the concept of a sequential, adaptive design for a single individual. To the best of our knowledge, this is the first work on learning the optimal individualized treatment rule in response to the current context for a single subject. A treatment rule for a patient is an individualized treatment strategy based on the history accrued, and context learned, up to the most current time point. A reward is measured on the patient at repetitive units, and optimality is meant in terms of optimization of the mean reward at a particular time $t$. We aim to learn an optimal, unknown choice of the controlled components of the design based on the data collected on the individual over time; with that, we adapt the randomization mechanism for future time-point experiments. Our results demonstrate that one can learn the optimal, context defined rule based on a single sample, and thereby adjust the design at any point $t$ with valid inference for the mean target parameter. 

This article provides two main contributions to the field of statistical precision medicine. First, we define a general class of averages of conditional (context-specific) causal parameters for the single unit time-series data. We define models for the probability distribution of the time-series that refrains from making unrealistic parametric assumptions, and aims to fully utilize the sequential randomization in the estimation procedure. Secondly, we present the study of the data-adaptive inference on the mean under the optimal treatment rule, where the target parameter adapts over time in response to the observed context of the individual. Our estimators are double robust and easier to estimate efficiently than previously proposed variations \citep{vanderLaan2018onlinets}. Finally, for inference, we rely on martingale Central Limit Theorem under a conditional variance stabilization condition and a maximal inequality for martingales with respect to an extension of the notion of bracketing entropy for martingale settings, initially proposed by \cite{vandeGeer2000}, which we refer to as \textit{sequential bracketing entropy}.

This article structure is as follows. In Section \ref{formulation} we formally present the general formulation of the statistical estimation problem, consisting of specifying the statistical model and notation, the target parameter defined as the average of context-specific target parameters, causal assumptions and identification results, and the corresponding efficient influence curve for the target parameter. In Section \ref{otr_ss} we discuss different strategies for estimating the optimal treatment rule and sampling strategies for assigning treatment at each time point. The following section, Section \ref{tmle}, introduces the Targeted Maximum Likelihood Estimator (TMLE), with Section \ref{theory} covering the theory behind the proposed estimator. In Section \ref{simulations} we present simulation results for different dependence settings. We conclude with a short discussion in Section \ref{conclusions}.

\section{Statistical Formulation of the Problem}
\label{formulation}

\subsection{Data and Likelihood}\label{sec:likelihood}

Let $O(t)$ be the observed data at time $t$, where we assume to follow a patient along time steps $t=1, \ldots, N$ such that $O^N \equiv (O(0), O(1), \ldots O(N)) = (O(t): t=0, \ldots, N)$. At each time step $t$, the experimenter assigns to the patient a binary treatment $A(t) \in \mathcal{A} := \{0,1\}$. 
%with our results trivially generalized to a vector $A(t)\in \{0,1\}^k$ of $k$ binary treatment nodes. 
We then observe, in this order, a post-treatment health outcome $Y(t) \in \mathcal{Y} \subset \mathbb{R}$, and then a post-outcome vector of time-varying covariates $W(t)$ lying in an Euclidean set $\mathcal{W}$. We suppose that larger values of $Y(t)$ reflect a better health outcome;
%, analogous in great generality to the reinforcement learning literature \citep{sutton1998}. 
without loss of generality, we also assume that $Y(t) \equiv (0,1)$, with rewards being bounded away from 0 and 1. The ordering of the nodes matters, as $W(t)$ is an important part of post-exposure history to be considered for the next record, $O(t+1)$. Finally, we note that $O(0) = (W(0))$, where $O(-1) = A(0) = Y(0) = \emptyset$; as such, $O(0)$ plays the role of baseline covariates for the collected time-series, based on which exposure $A(1)$ might be allocated. 

We denote $O(t) := (A(t), Y(t), W(t))$ the observed data collected on the patient at time step $t$, with $\mathcal{O} := \mathcal{A} \times \mathcal{Y} \times \mathcal{W}$ as the domain of the observation $O(t)$. We note that $O(t)$ has a fixed dimension in time $t$, and is an element of an Euclidean set ${\cal O}$. Our data set is the time-indexed sequence $O^N \in \mathcal{O}^N$, or \textit{time-series}, of the successive observations collected on a single patient. For any $t$, we let $\bar{O}(t) := (O(1),\ldots,O(t))$ denote the observed history of the patient up until time $t$. Unlike in more traditional statistical settings, the data points $O(1),\ldots,O(N)$ are not independent draws from the same law: here they form a dependent sequence, which is a single draw of a distribution over $\mathcal{O}^N$. In that sense, our data reduces to a single sample. 

We let $O^N \sim P_0^N$, where $P_0^N$ denotes the true probability distribution  of $O^N$. The subscript ``0'' stands for the ``truth'' throughout the rest of the article, denoting the true, unknown features of the distribution of the data. 
Realizations of a random variable $O^N$ are denoted with lower case letters, $o^N$. 
%For consistency, we also define $Y^N = (Y(t): t=1, \ldots, N)$, $A^N = (A(t): t=1, \ldots, N)$ and $W^N = (W(t): t=0, \ldots, N)$. 
%Coarsened series from time point $j$ to $N$ is written as $O^{j:N} = (O(t): t=j, \ldots, N)$, whereas the past observations from time $t$ until the first time-point are denoted as $\overline{o}(t) = (o(0), \ldots, o(t))$. 
We suppose that $P_0^N$ admits a density $p^N_0$ w.r.t. a dominating measure $\mu$ over $\mathcal{O}^N$ that can be written as the product measure $\mu = \times_{t=1}^N (\mu_A \times \mu_Y \times \mu_W)$, with $\mu_A$, $\mu_Y$, and $\mu_W$ measures over $\mathcal{A}$, $\mathcal{Y}$, and $\mathcal{W}$. From the chain rule, the likelihood under the true data distribution $P^N_0$ of a realization $\bar{o}^N$ of $\bar{O}^N$ can be factorized according to the time ordering of observation nodes as:
\begin{align}\label{ll} 
p_0^N(o^N) = &\prod_{t=1}^N p_{0, a(t)}(a(t) \mid \overline{o}(t-1)) \times \prod_{t=1}^N p_{0, y(t)}(y(t) \mid \overline{o}(t-1), a(t)) \\ 
& \times \prod_{t=0}^N p_{0, w(t)}(w(t) \mid \overline{o}(t-1), a(t), y(t)),
\end{align}
where $a(t) \mapsto p_{0,a(t)}(a(t) \mid \bar{o}(t-1)$, $y(t) \mapsto p_{0,y(t)}(y(t) \mid \bar{o}(t-1), a(t))$, and $w(t) \mapsto p_{0,w(t)}(w(t) \mid \bar{o}(t-1), a(t), y(t))$ are conditional densities w.r.t. the dominating measures $\mu_A$, $\mu_Y$, and $\mu_W$. 

%Let $\mu_W$, $\mu_A$ and $\mu_Y$ be measures on the outcome spaces $\mathcal{W}$, $\mathcal{A}$ and $\mathcal{Y}$ for $(W(t),A(t),Y(t))$, for any $t$, equipped with the appropriate $\sigma$-fields. In addition, we define $\mu = \mu_W \times \mu_A \times \mu_Y$ as the measure on $\mathcal{O}$ equipped with the product of the three $\sigma$-fields. The unknown, true likelihood of $O^N$ with respect to $\mu$ is given by the following factorization of the probability density according to the time ordering:
%\begin{align}\label{ll} 
%p_0^N(o^N) = &\prod_{t=1}^N p_{0, a(t)}(a(t) \mid \overline{o}(t-1)) \prod_{t=1}^N p_{0, y(t)}(y(t) \mid \overline{o}(t-1), a(t)) \\ &\prod_{t=0}^N p_{0, w(t)}(w(t) \mid \overline{o}(t-1), a(t), y(t)) \nonumber
%\end{align} 

%\noindent where $w \rightarrow p_{0, w(t)}( \ \cdot \ | \overline{o}(t-1),a(t),y(t))$ is the conditional density w.r.t dominating
%measure $\mu_W$ of a true, unknown law $\mathcal{W}$, and $Y \rightarrow
%p_{0,y(t)}( \ \cdot \ | \overline{o}(t-1),a(t))$ is the true conditional probability density with respect to dominating measure $\mu_{Y}$. 
%We define a vector $g^N \equiv (g_1, \ldots g_N)$ as the ordered vector of the first $N$ stochastic rules, with $p_{a(t)}( \ \cdot \ \mid \overline{o}(t-1)) = g_t(\ \cdot \ \mid \overline{o}(t-1))$ denoting the known conditional probability of $t$-specific treatment given the history. 

\subsection{Statistical Model}\label{sec:statmodel}

Since $O^N$ represents a single time-series, a dependent process, we observe only a single draw from $P_0^N$. As a result, we are unable to estimate any part of $P_0^N$ without additional assumptions. In particular, we assume that the conditional distribution of $O(t)$ given $\bar{O}(t-1)$, $P_{O(t)|\bar{O}(t-1)}$, depends on $\bar{O}(t-1)$ through a summary measure $C_o(t)=C_o(\bar{O}(t-1))\in {\cal C}$ of fixed dimension; each $C_o(t)$ might contain $t$-specific summary of previous measurements of context, or is of a particular Markov order. For later notational convenience, we denote this conditional distribution $P_{O(t)\mid \bar{O}(t-1)}$ with $P_{C_o(t)}$. Then, the density $p_{C_o(t)}$ of $P_{C_o(t)}$ with respect to a dominating measure $\mu_{C_o(t)}$ is a conditional density $(o,C_o)\rightarrow p_{C_o(t)}(o\mid C_o)$ so that for each value of $C_o(t)$, $\int p_{C_o(t)}(o\mid C_o(t))d\mu_{C_o(t)}(o)=1$. We extend this notion to all parts of the likelihood as described in subsection (\ref{ll}), defining $q_{y(t)}$ as the density for node $Y(t)$ conditional on a fixed dimensional summary $C_y(t)$, with $C_w(t)$ and $C_a(t)$ corresponding to fixed dimensional summaries for $q_{w(t)} = p_{0,w(t)}(w(t) \mid C_w(t))$ and $g_t = p_{0,a(t)}(a(t) \mid C_a(t))$, respectively. 

Additionally, we assume that $p_{C_o(t)}$ is parameterized by a common (in time $t$) function $\theta\in \Theta$, with inputs $(c,o)\rightarrow \theta(c,o)$. The conditional distribution $p_{C_o(t)}$ depends on $\theta$ only through $\theta(C_o(t),\cdot)$
%observed $C_o(t)$, with $\theta \equiv \theta(C_o(t),\cdot)$.
We write $p_{C_o(t)} = p_{\theta, C_o(t)}$ interchangeably.
%as $\theta$ always implies dependence on $C_o(t)$.
Let $q_y$ be the common conditional density of $Y(t)$, given $(A(t),C_o(t))$; we make no such assumption on $q_{w(t)}$. Additionally, we make no conditional stationarity assumptions on $g_t$ if randomization probabilities are known, as is the case for an adaptive sequential trial. We define $\bar{Q}(C_o(t), A(t))=E_{P_{C_o(t)}}(Y(t)\mid C_o(t),A(t))$ to be the conditional mean of $Y(t)$ given $C_o(t)$ and $A(t)$. As such, we have that $\bar{Q}(C_y(t)) = \bar{Q}(C_o(t), A(t))=\int y q_y(y\mid C_o(t),A(t))d\mu_y(o)$, and $\bar{Q}$ is a common function across time $t$; we put no restrictions on $\bar{Q}$. We suppress dependence of the conditional density $q_{w(t)}$ in future reference, as this factor plays no role in estimation. In particular, $q_{w(t)}$ does not affect the efficient influence curve of the target parameter, allowing us to act as if $q_{w(t)}$ is known.
Finally, we define $\theta=(g,\bar{Q})$ and let $\Theta={\cal G}\times \bar{\mathcal{Q}}$ be the cartesian product of the two nonparametric parameter spaces for $g$ and $\bar{Q}$. 
%We choose to define $\theta$ as $\theta=(g,\bar{Q})$ in order to include more general scenarios ($g$ not randomized) in further definitions. 

Let $p_{\theta,C_o(t)}$ and $p_{\theta}^N$ be the density for $O(t)$ given $C_o(t)$ and $O^N$, implied by $\theta$. This defines a statistical model ${\cal M}^N=\{P^N_{\theta}:\theta\}$ where $P^N_{\theta}$ is the probability measure for the time-series implied by 
$p_{\theta, C_o(t)}$. Additionally, we define a statistical model of distributions of $O(t)$ at time $t$, conditional on realized summary $C_o(t)$. In particular, let ${\cal M}(C_o(t))=\{P_{\theta,C_o(t)}:\theta\}$ be the model for $P_{C_o(t)}$ for a given $C_o(t)$ implied by ${\cal M}^N$. Note that, by setup, both ${\cal M}^N$ and  ${\cal M}(C_o(t))$ contain their truth $P_0$ and $P_{O(t)\mid C_o(t)}$, respectively. Similarly to the likelihood expression in sub section (\ref{ll}), we can factorize the likelihood under the above defined statistical model according to time ordering as:
\begin{align}\label{l2} 
p_{\theta}(o^N) = \prod_{t=1}^N g_t(a(t) \mid C_{a}(t)) \prod_{t=1}^N {q}_y(y(t) \mid C_{y}(t)) \prod_{t=0}^N q_{w(t)}(w(t) \mid C_{w}(t)).
\end{align} 

\subsection{Causal Target Parameter and Identification}\label{sec:causaltarget}

\subsubsection{Non-parametric structural equation model and causal target parameter}

%\paragraph{Non-parametric structural equation model.} 
By specifying a non-parametric structural equations model (NPSEM; equivalently, structural causal model), we assume that each component of the observed time-specific data structure is a function of an observed, fixed-dimensional history and an unmeasured exogenous error term (\cite{pearl2009}). We encode the time-ordering of the variables using the following NPSEM:
\begin{align}\label{eq:npsem}
W(0) &= f_{W(0)}(U_{W}(0)), \\ 
A(t) &= f_{A(t)}(C_{A}(t), U_{A}(t)), \ \ \ t=1, \ldots, N, \\ 
Y(t) &= f_{Y(t)}(C_{Y}(t), U_{Y}(t)), \ \ \ t=1, \ldots, N, \\ 
W(t) &= f_{W(t)}(C_{W}(t), U_{W}(t)), \ \ t=1, \ldots, N,
\end{align}
where $(f_{A}(t): t=1, \ldots, N)$, $(f_{Y}(t): t=1, \ldots, N)$ and $(f_{W}(t): t=0, \ldots, N)$ are unspecified, deterministic functions and $U = (U_{W}(0), \ldots, U_{W}(N),  U_{A}(1), \ldots,  U_{A}(N), U_{Y}(1), \ldots, U_{Y}(N))$ is a vector of exogenous errors.

We denote $\mathcal{M}^F$ the set of all probability distributions $P^F$ over the domain of $(O,U)$ that are compatible with the NPSEM defined above. Let $P^F_0$ be the true probability distribution of $(O,U)$, which we assume to belong to $\mathcal{M}^F$; we denote $\mathcal{M}^F$ as the \textit{causal model}.
%The causal model, $\mathcal{M}^F$, reflects a model for the distribution of $(O,U)$; it provides a parameterization of the distribution of the observed data structure in terms of full-data distribution of $(O,U)$, including both measured and unmeasured variables, modeled by the system of structural equations.
The causal model  $\mathcal{M}^F$ encodes all the knowledge about the data-generating process, and implies a model for the distribution of the counterfactual random variables; as such, causal effects are defined in terms of hypothetical interventions on the NPSEM. 

Consider a treatment rule $C_o(t) \rightarrow d(C_o(t)) \in \{0,1\}$, that maps the observed, fixed dimensional history $C_o(t)$ into a treatment decision for $A(t)$.
%with certain probability. 
We introduce a counterfactual random variable $O^{N,d}$, defined by substituting the equation for node $A$ at time $t$ in the NPSEM with the intervention $d$:
\begin{align}\label{eq:counterfactual_NPSEM}
    W^d(0) &= f_{W(0)}(U_{W}(0)) \\
    %A^d(t) &\sim g^*(\cdot \mid C_{A}(t)), \ \ \ t=1, \ldots, N \\
    A^d(t) &= d(C_{A}(t)), \ \ \ t=1, \ldots, N \\
    Y^d(t) &= f_{Y(t)}(C_{Y}(t), U_{Y}(t)), \ \ \ t=1, \ldots, N \\
    W^d(t) &= f_{W(t)}(C_{W}(t), U_{W}(t)), \ \ t=1, \ldots, N,
\end{align}
%where $g^*$ is the conditional distribution of treatment assignment given the past summary vector, under stochastic treatment rule $d$. 
We gather all of the nodes of the above modified NPSEM in the random vector $O^{N,d} := (O^d(t) : t =1,\ldots,N)$, where $O^d(t) := (A^d(t), Y^d(t), W^d(t))$. The random vector $O^{N,d}$ represents the counterfactual time-series, or counterfactual trajectory the subject of interest would have had, had each treatment assignment $A(t)$, for $t=1,\ldots,N$, had been carried out following the 
%stochastic 
treatment rule $d$.

% The post-intervention distribution in is defined as the distribution $O^{N,d}$ would have had, had the treatment been assigned according to the $t^{\text{th}}$-time-point specific rule, $d(C_o(t))$ across all time points $t$.

%\paragraph{Causal target parameters} 

%The main way in which the present work differs from other articles from the time-series causal inference literature is that we consider time- and context-specific causal parameters, which as we show below, leads to a host of desirable statistical properties. 

We now formally define time-series causal parameters. First, we introduce a time- and context-specific causal model.
Let $\mathcal{M}^F_t(c)$ be the set of conditional probability distributions $P^F_{c}$ over the domain of $(O(t), U_A(t), U_Y(t), U_W(t))$ compatible with the non-parametric structural equation model \eqref{eq:npsem} imposing that $C_A(t) = C_o(t) = c$:
\begin{align}
    A_c(t) =&  f_{A(t)}(c_o, U_A(t)) \\
    Y_c(t) =& f_{Y(t)}(c_A(c_o, A(t)), U_Y(t)) \\
    W_c(t) =& f_{W(t)}(c_W(c_o, A(t), Y(t)), U_W(t)).
\end{align}
Let $O^d_{c}(t)$ be the counterfactual observation at time $t$, obtained by substituting the $A(t)$ equation in the above set of equations with the deterministic intervention $d$:
\begin{align}
    %A^d_c(t) \sim& g^*(\cdot \mid c) \\
    A^d_c(t) =& d(c) \\
    Y^d_c(t) =& f_{Y(t)}(c_A(c_o, A(t)), U_Y(t)) \\
    W^d_c(t) =& f_{W(t)}(c_W(c_o, A(t), Y(t)), U_W(t)).
\end{align}
We define our causal parameter of interest as
\begin{align}
    \Psi^{F,d}_{C_o(t)}(P^F_{C_o(t)}) := E[Y^d_{C_o(t)}],
\end{align}
which is the expectation of the counterfactual random variable $Y^d$, generated by the above modified NPSEM. It corresponds to starting at $c = C_o(t)$, the current context, and assigning treatment following $d$.
%$g^*$ - the conditional treatment distribution corresponding to stochastic intervention $d$.
Our causal target parameter is the mean outcome we would have obtained after one time-step, if, starting at time $t$ from the observed context $C_o(t)$, we had carried out intervention $d$.

\subsubsection{Identification of the causal target and defining the statistical target}

Once we have defined our causal target parameter, the natural question that arises is how to identify it from the observed data distribution. We can identify the distribution of the $d$-specific time series $O^{N,d}$, and also of the $(d, C_o(t))$-specific observation $O_{C_o(t)}^d$, from the observed data via the G-computation formula - under the sequential randomization and positivity assumptions, which we state below.

\begin{assumption}[Sequential randomization]\label{assumption:sequential_randomization}
For every $t$, $Y^d(t) \ \indep \ A(t) \mid C_o(t)$ (and $Y^d_{C_o(t)}(t) \ \indep \ A(t) \mid C_o(t)$).
\end{assumption}

\begin{assumption}[Positivity]\label{assumption:positivity}
It holds that under the treatment mechanism $g_{0,t}$, each treatment value $a \in \{0,1\} $ has a positive probability of being assigned, under every possible treatment history:
\begin{align}
    g_{0,t}(a \mid c) > 0, \text{ for every } t \geq 1,\ a \in \{0,1\} \text{ and every } c \in \mathcal{C} \text{ such that } P_0[C_o(t) = c] > 0.
\end{align}
\end{assumption}

Note that under the setting of the present article, as we suppose that $A(t)$ is assigned at random conditional on $C_o(t)$ by the experimenter, assumption 1 concerning the sequential randomization automatically holds. Under identification assumptions 1 and 2, we can write our causal parameter $\Psi^{F,d}_{C_o(t)}(P^F_0)$ as a feature of the data-generating distribution:
%\begin{align}
%    \Psi^{F,d}_{C_o(t)}(P^F_0) = \Psi^d_{C_o(t)}(P_0) := E_{P_0} \left[ Y(t) \frac{g^*(A(t) \mid C_o(t))}{g_{0,t}(A(t) \mid C_o(t))} \mid C_o(t) \right],
%\end{align}
%which for a deterministic $d$ can be written as 
\begin{align}
    \Psi^{F,d}_{C_o(t)}(P^F_0) = \Psi_{C_o(t)}^d(P_0) := E_{P_0} \left[ Y(t) \mid A(t) = d(C_o(t)), C_o(t) \right].
\end{align} 
\noindent
Observe that in the above definition $\Psi^d_{C_0(t)}(P_0)$, shortly defined as $\Psi_{C_0(t)}(P_0)$,  depends on $P_0$ only through the true conditional distribution of $O(t)$ given $C_o(t)$. For every $P$, we remind that $P_{C_o(t)}$ denotes the distribution of $O(t)$ given $C_o(t)$, and let $\mathcal{M}(C_o(t))$ be the set of such distributions corresponding to $P \in \mathcal{M}$. At each time-point $t$, given a $C_o(t)$, we define a target parameter $\Psi_{C_o(t)} : \mathcal{M}(C_o(t)) \rightarrow \mathbb{R}$ that is pathwise differentiable with canonical gradient $D^*_{C_o(t)}(P_{C_o(t)})(o)$ at $P_{C_o(t)}$ in $\mathcal{M}(C_o(t))$. As described in Section \ref{sec:statmodel}, we have that $\Psi_{C_o(t)}(P_{C_o(t)}) =  \Psi_{C_o(t)}(\theta)$, where $\Psi_{C_o(t)}(\theta)$ depends on $\theta$ only though its section $\theta(C_o(t), \cdot)$. We denote the collection of $C_o(t)$-specific canonical gradients as $(c,o) \rightarrow D^*(P_{C_o(t)})(c,o)$, so that we can write them uniformly as a function of the observed components; with that, we have that $D^*_{C_o(t)}(P_{C_o(t)})(o) = D^*_{C_o(t)}(\theta)(o) = D^*(\theta)(c_o(t), o)$. As is custom for canonical gradients, for a given $C_o(t)$, $D^*(\theta)$ is a function of the observed data with conditional mean zero with respect to $P_{C_o(t)}$.

Finally, we propose a class of statistical target parameters $\bar{\Psi}(\theta)$ defined as the average over time of $C_o(t)$-specific counterfactual means under the treatment rule. In particular, the target parameter on ${\cal  M}^N$, $\Psi^N:{\cal M}^N \rightarrow \mathbb{R}$ of the data distribution $P^N \in {\cal M}^N$ is defined as:
\begin{equation}\label{tp2}
    \bar{\Psi}(\theta)=\frac{1}{N} \sum_{t=1}^N \Psi_{C_o(t)}(\theta).
\end{equation}
The statistical target parameter $\bar{\Psi}(\theta)$ is data-dependent, as it is defined as an average over time of parameters of the conditional distribution of $O(t)$ given the observed realization of $C_o(t)$; as such, it depends on $(C_o(1), \ldots, C_o(N))$. In practice, $\bar{\Psi}(\theta)$ is an average of the means under optimal treatment decisions over all observed contexts over time. As an average of $C_o(t)$-specific causal effects with a double robust efficient influence curve $D^*_{C_o(t)}(\theta)(o)$, it follows we can estimate $\bar{\Psi}(\theta)$ in a double robust manner as well, as we further emphasize in the following section. 

%In an adaptive sequential design, the treatment is sequentially randomized according to known treatment probabilities at each time $t$, hence the sequential randomization assumption is automatically satisfied. We also note that assumptions (1) and (2) hold by construction in a micro-randomized trail, with (2) also holding automatically for stochastic treatment rules \citep{klasnja2015,chambaz2017}. With that, under assumptions (1)-(2), the average over time of $C_o(t)$-specific parameters can be expressed in terms of the observed data. 

\subsubsection{Canonical gradient and first order expansion of the target parameter}

%The key theoretical result of this article is the derivation of the  canonical gradient of our target parameter that admits a first order expansion with a double-robust second order term. 

In the following theorem we provide the canonical gradient of our target parameter that admits a first order expansion with a double-robust second order term.  We pursue the discussion on marginal parameters in more detail in subsection \ref{subsection:comparison_marg_params} in the appendix.

\begin{theorem}[Canonical gradient and first order expansion]\label{thm:can_gdt_and_1st_order_exp}
Under the strong positivity assumption (assumption \ref{assumption:strong_positivity} in subsection \ref{M1}), the target parameter mapping $\Psi_{C_0(t)} : \mathcal{M}(C_o(t)) \to \mathbb{R}$ is pathwise differentiable w.r.t. $\mathcal{M}(C_o(t))$, with a canonical gradient w.r.t. $\mathcal{M}(C_o(t))$ given by
\begin{align}
    D^*_{C_o(t)}(\bar{Q}, g)(o) = \frac{g^*(a \mid C_o(t))}{g(a \mid C_o(t))} \left(y - \bar{Q}(a, C_o(t))\right).
\end{align}
Furthermore $\Psi_{C_o(t)}(\bar{Q})$ admits the following first order expansion:
\begin{align}
    \Psi_{C_o(t)}(\bar{Q}) - \Psi_{C_o(t)}(\bar{Q}_0) = - P_{0,C_o(t)} D^*_{C_o(t)}(\bar{Q}, g) + R(\bar{Q},\bar{Q}_0, g, g_{0,t}),
\end{align}
where $R$ is a second order remainder that is doubly-robust, with
$R(\bar{Q},\bar{Q}_0, g, g_{0,t})=0$ if either $\bar{Q} = \bar{Q}_0$ or $g = g_{0,t}$.
\end{theorem}

%Previous works on statistical parameters defined over a single time series model \citep{vanderLaan2018onlinets, kallus2019efficiently} consider what we refer to as \textit{marginal} parameters. Unlike the conditional parameters we consider here, the efficient influence function of marginal parameters is not double-robust in the usual sense; that is, robust w.r.t. a pair of variation independent nuisance parameters. More importantly, knowing or consistently estimating the treatment mechanism does not guarantee consistency of the causal effect for parameters described by \cite{vanderLaan2018onlinets} and \cite{ kallus2019efficiently}. We pursue the discussion on marginal parameters in more detail in subsection \ref{subsection:comparison_marg_params} in the appendix.

\subsubsection{Optimal rule}

Now that we have identified the context-specific counterfactual outcome under $d$ as a parameter of the observed data distribution $P_0^N$, we can identify the optimal treatment rule. The  optimal treatment rule is a priori a causal object defined as a function of $P^F_0$, and a  parameter of the observed data generating distribution $P_0^N$. Under the identification assumptions, we can identify the optimal rule from the observed data distribution as follows. Fix arbitrarily $\bar{Q} \in \bar{\mathcal{Q}}$. To alleviate notation, we further introduce the blip function defined as:

\begin{equation}\label{blip}
    B(C_o(t)) \equiv \bar{Q}(C_o(t), A(t)=1) - \bar{Q}(C_o(t), A(t)=0).
\end{equation}
Intuitively, if $B(C_o(t)) > 0$, assigning treatment $A(t) = 1$ is more beneficial (in terms of optimizing $Y(t)$) than $A(t) = 0$ for time point $t$ under the current context $C_o(t)$. If $B(C_o(t)) < 0$, we can optimize the $t$-specific outcome by assigning the subject treatment $A(t) = 0$ instead. The true optimal rule for the purpose of optimizing the mean of the next (short-term) outcome $Y(t)$, for binary treatment, is then given by:
\begin{equation}\label{eq:optdef}
    d_{0}(C_o(t)) \equiv \mathbb{I}(B(C_o(t)) > 0).
\end{equation}
As defined in Equation \eqref{eq:optdef}, $d_{0}(C_o(t))$ is a typical treatment rule that maps observed fixed dimensional summary deterministically into one treatment; a stochastic treatment rule does so randomly \citep{luedtke2016,luedtke2016resource,chambaz2017}. 
%We discuss potential disadvantages of deterministic treatment rules in the following section, as well as several strategies for learning, and approximating, optimal stochastic rules.

%To put it in more context, consider the $C_o(t)$-specific conditional counterfactual mean under the optimal stochastic rule $d(C_o(t))$. 
%In the following, we introduce two targets for the adaptive sequential design. We focus on exposures with proximal (in time) effects, dependent on the current context and history of the subject, as often dictated by mHealth applications (\cite{heron2010}). 

%As motivation, consider Sense2Stop - a micro-randomized trial aiming to investigate whether sensor-based assessment of stress, and appropriate interventions to relax, can help chronic smokers from relapsing \citep{sense2stop}. In practice, we may choose to follow a single individual enrolled in the Sense2Stop MRT, and at each time $t$ allocate the stress-management intervention based on some summary of the past. For example, we may have learned over time that the individual is more stressed on certain days, or in proximity to certain locations. In addition, the dense time-series data may have elucidated that the immediate urge to smoke is lessened for this participant if the intervention is executed outside the working hours, no matter what the stress level was. 

%%%%%%%%%%%%%%%%%%%%%%%%%%%%%%%%%%%%%%%%%%%%%%%%%%
% Learning the Optimal Rule
%%%%%%%%%%%%%%%%%%%%%%%%%%%%%%%%%%%%%%%%%%%%%%%%%%

\section{Optimal Treatment Rule and the Sampling Scheme}\label{otr_ss}

In an adaptive sequential trial, the process of generating $A(t)$ is controlled by the experimenter. As such, one can simultaneously learn and start assigning treatment according to the best current estimate of the optimal treatment rule, with varying exploration-exploitation objectives. In this section we describe different strategies for estimating the optimal treatment rule, as well as propose different sampling schemes for assigning treatment.  
%define $g \equiv (g_1, \ldots, g_n)$ --- a targeted time-specific scheme that "targets" the optimal individualized treatment. 

\subsection{Estimating the Optimal Treatment Rule}

%Perhaps a bit confusing to parametrize with \theta here?
\subsubsection{Parametric working model}
First, we consider estimating the optimal treatment rule based on a parametric working model. As described previously, consider a treatment rule $C_o(t)\rightarrow d(C_o(t))\in \{0,1\}$ that maps the history $C_o(t)$ into a treatment decision for $A(t)$. We define a parametric working model for $q_y$ indexed by parameter $\theta$ such that $\{q_{y,\theta}: \theta\}$. Notice that under the specified working model, we have that:
\begin{equation*}
\bar{Q}_{\theta}(C_o(t),a) = E(Y(t)\mid C_o(t),A(t)=a)=\int y q_{y,\theta}(y \mid C_o(t),a) d\mu_y(y).
\end{equation*}
\noindent
The true conditional treatment effect, $B_0(C_o(t))$, can then be expressed as
\begin{equation*}
B_{\theta}(C_o(t))=\bar{Q}_{\theta}(C_o(t),1)-\bar{Q}_{\theta}(C_o(t),0)
\end{equation*}
under the parametric working model. Recall that the optimal treatment rule for $A(t)$ for the purpose of maximizing $Y(t)$ is given by: 
\begin{equation*}
    d_0(C_o(t))=\mathbb{I}(B_0(C_o(t))>0).
\end{equation*}
Under the parametric working model, we note that the optimal treatment rule can be represented as: 
\begin{equation*}
    d_{\theta}(C_o(t))=\mathbb{I}(B_{\theta}(C_o(t))>0).
\end{equation*}

Let $\theta_{t-1}$ to be the maximum likelihood estimate of the true $\theta_0$ based on the most current history, $\bar{O}(t-1)$, and according to the working model $q_{y,\theta}$. We could define the fixed dimensional history $C_o(t)$ such that for each time point $t$, $\theta_{t-1}$ is included in the relevant history $C_o(t)$ for $O(t)$. The current estimate of the rule is then defined as:
\begin{equation*}
    d_{\theta_{t-1}}(C_o(t))=\mathbb{I}(B_{\theta_{t-1}}(C_o(t))>0).
\end{equation*}
\noindent
If the parametric model is very flexible, $B_{\theta_{t-1}}$ might be a good approximation of the true conditional treatment effect $B_0(C_o(t))$. In that case, $d_{\theta_{t-1}}(C_o(t))$ is a good approximation of the optimal rule $d_0(C_o(t))$. Nevertheless, we argue that $\theta_{t-1}$ will converge to $\theta_0$ defined by a Kullback-Leibler projection of the true $q_{y,0}$ onto the working model $\{q_{y,\theta}:\theta\}$. Consequently, the rule $d_{\theta_{t-1}}(C_o(t))$ will converge to a fixed  $\mathbb{I}(B_{0}(C_o(t))>0)$  as $t$ converges to infinity.

\subsubsection{Machine Learning}

Instead of considering a parametric working model, we explore estimation of the optimal treatment rule based on more flexible, possibly nonparametric approaches drawn from the machine learning literature. As in the previous subsection, we define $B_{\bar{Q}_{t-1}}(C_o(t))$ to be an estimator of the true blip function, $B_0(C_o(t))$, based on the most recent observations up to time $t$, $\bar{O}(t-1)$. In particular, we consider estimators studied in previous work, including online Super-Learner of $\bar{Q}_0$ which provides convenient computational and statistical properties for dense time-series data described elsewhere \citep{lendle2014,benkeser2018}. Additionally, we might consider ensemble machine learning methods that target $B_0$ directly \citep{luedtke2016}. As mentioned in the previous section, we can view $B_{\bar{Q}_{t-1}}(C_o(t))$ as just another univariate covariate extracted from the past, and include it in our definition of $C_o(t)$.  If $B_{\bar{Q}_{t-1}}$ is consistent for $B_0$, then the rule $d_{\bar{Q}_{t-1}}(C_o(t))$ based on $B_{\bar{Q}_{t-1}}$ will converge to the optimal rule $\mathbb{I}(B_0(C_o(t))>0)$, as shown in previous work \citep{luedtke2016, chambaz2017}.

\subsection{Defining the Sampling Scheme}

In the following section, we describe two sampling schemes that define $g^N = \{g_t: t = 1, \cdots, N\}$ precisely. Both rely on estimating parts of the likelihood based on the time-points collected so far for the single subject studied. The $t$-dependent current estimate of $\bar{Q}_0$ and $B_0$ are then further utilized to assign the next treatment, collect the next  corresponding block of data, and estimate the target parameter of interest. Following the empirical process literature, we define $P_N f$ to be the empirical average of function $f$, and $Pf = \mathbb{E}_Pf(O)$.

\subsubsection{Stochastic Optimal Treatment Rules}\label{sotr}

In the following, we follow closely the argument given in \cite{chambaz2017} in order to define $g^N$. Let $\bar{Q}_{t-1}$ denote the time $t$ estimate of $\bar{Q}_0$ based on the time-series points collected so far, $\bar{O}(t-1)$. For a small number of samples, $d_{\bar{Q}_{t-1}}(C_o(t))$ might not be a good estimate of $d_0(C_o(t))$. As such, assigning the current conditional
probability of treatment deterministically based on the estimated rule could be
ill-advised. In addition, without exploration (enforced via a deterministic rule), we cannot guarantee consistency of the optimal rule estimator. 

In light of that, we define $\{c_t\}_{t \geq 1}$ and $\{e_t\}_{t \geq 1}$ as user-defined, non-increasing sequences such that $c_1 \leq
\frac{1}{2}$, $\lim_{t} c_t \equiv c_{\infty} > 0$ and $\lim_{t} e_{t} \equiv
e_{\infty} > 0$. More specifically, we let $\{e_t\}_{t \geq 1}$ define the
level of random perturbation around the current estimate $d_{\bar{Q}_{t-1}}(C_o(t))$ of the optimal rule. We define $\{c_t\}_{t \geq 1}$ as the probability of failure, so choosing $c_1 = \cdots = c_t = 0.5$ would yield a balanced stochastic treatment rule. In particular, we define a design that ensures that, under any context and with a positive probability $e_t$, we pick the treatment uniformly at random. This positive probability $e_t$ is what is often referred to as the \textit{exploration rate} in the bandit and reinforcement learning literature \citep{sutton1998}.
For every $t \geq 1$, we could have the following function $G_t$ over $[-1,1]$ as defined in \cite{chambaz2017}:
\begin{equation*}
G_t(x) = c_t  \mathbb{I}[x \leq -e_t] + (1-c_t) \mathbb{I}[x \geq e_t] +
\left(-\frac{1/2-c_t}{2 e_t^3}x^3 + \frac{1/2-c_t}{2 e_t/3}x +
  \frac{1}{2}\right) \mathbb{I}[- e_t \leq x \leq e_t],
\end{equation*}
where $G_t(x)$ is used to derive a stochastic treatment rule
from an estimated blip function, such that
\begin{equation*}
    g_t(1 \mid C_o(t)) = G_t(B_{\bar{Q}_{t-1}}(C_o(t))).
\end{equation*}

Note that $G_t$ is a smooth approximation to $x \rightarrow \mathbb{I}[x \geq 0]$
bounded away from 0 and 1, mimicking the optimal treatment rule as an indicator
of the true blip function. With that in mind, any other non-decreasing $k_n$-Lipschitz function with $F_t(x) = c_t$ for $x \leq -e_t$ and $F_t(x) = 1 - c_t$ for $x \geq e_t$ would approximate the optimal treatment rule as well. The definitions of $G_t$ and $g_t$ prompt the following lemma, which illustrates the ability of the sampling scheme to learn form the collected data, while still exploring:

\vspace{3mm} \noindent 
\begin{lemma}\label{eeto}
%\textbf{Exploration/Exploitation Trade Off} \\
\textit{Let t $\geq$ 1. Then we have that: 
\begin{align*}
    \inf_{c_o(t)} g_t(d(c_o(t)) \mid c_o(t)) &\geq \frac{1}{2} \\
    \inf_{c_o(t)} g_t(1-d(c_o(t)) \mid c_o(t)) &\geq c_t. 
\end{align*}}
\end{lemma}
\noindent
Note that under Lemma \ref{eeto}, the positivity assumption needed for the identification result is met. Finally, we reiterate that the stochastic treatment rule $g_t(1 \mid C_o(t))$ approximates $d(C_o(t))$ in the following sense:
\begin{equation*}
    |g_t(1 \mid C_o(t)) - d(C_o(t))| \leq c_{\infty} \mathbb{I}[|B(C_o(t)) \geq e_{\infty}|] + \frac{1}{2} \mathbb{I}[|B(C_o(t)) < e_{\infty}|].
\end{equation*}
If $c_{\infty}$ and $e_{\infty}$ are small and $|B(C_o(t)) \geq e_{\infty}|$, then drawing treatment assignment from a smooth approximation of $d(C_o(t))$ is not much different than $d(C_o(t))$, with little impact on the mean value of the reward. 

\subsubsection{Target sequential sampling with Highly Adaptive Lasso}

Alternatively, one could allocate randomization probabilities based on the tails of an estimate of the blip function, $B(C_o(t))$. In particular, we present a sampling scheme that utilizes the Highly Adaptive Lasso (HAL) estimator for obtaining the bounds around the estimate of the true blip function. The Highly Adaptove Lasso is a nonparametric regression estimator that does not rely on local smoothness assumptions \citep{benkeser2016,vanderlaan2017hal}. Briefly, for the class of functions that are right-hand continuous with left-hand limits and a finite variation norm, HAL is an MLE which can be computed based on $L_1$-penalized regression. As such, it is similar to standard lasso regression function in its implementation, except that the relationship between the predictors and the outcome is described by data-dependent basis functions instead of a parametric model. For a thorough description of the Highly Adaptive Lasso estimator, we refer the reader to \cite{benkeser2016} and \cite{vanderlaan2017hal}; we provide more details on the Highly Adaptive Lasso in the appendix subsection \ref{halmore}.

We propose to use HAL to estimate $B_0(C_o)$, which implies
an estimator for the optimal rule $d_0(C_o)=\mathbb{I}(B_0(C_o)>0)$. We define a quadratic loss function as follows: 
\begin{equation*}
   L_B(\theta)(o, C_o) =(D_1(\theta)(O) - B(C_o))^2, 
\end{equation*}
which is indexed by $\theta = (g,\bar{Q})$ required to evaluate $D_1(\theta)(O)$.
%\begin{equation*}
% \frac{2(A-1)}{g(A|C_o)}(Y-\bar{Q}(C_o,A))+ \bar{Q}(C_o,1) + \bar{Q}(C_o,0).
%\end{equation*}
This influence function has the property that $E_0(D_1(\theta)|C_o) = B_0(C_o)$ if either $\bar{Q} = \bar{Q}_0$ or $g = g_0$, under positivity. As such, $L_{B}(\theta_0)$ is a double robust and efficient loss function for the true risk in the sense that $P_n L_{B}(\theta)$ is a double robust locally efficient estimator of the true risk under regularity conditions. As a double robust and efficient loss, the true risk of the loss function $L_B(\theta)(o,C_o)$ equals $P_0 (B_0-B)^2 (C_o)$ up until a constant if either $D_1(\theta) = D_1(\bar{Q}_0,g)$ or $D_1(\theta) = D_1(\bar{Q}, g_0)$.

Let $E(D_1(\theta) | C_o) = \psi^{\text{blip}}$, with $\psi^{\text{blip}} \in D[0, \tau]$, the Banach space of $d$-variate cadlag functions. Define $C_{o,s} = \{C_{o,j} : j \in s\}$ for a given subset $s \subset \{1, ..., d\}$. For $\psi^{\text{blip}} \in D[0,\tau]$, we define the $s^{\text{th}}$ section of $\psi^{\text{blip}}$ as $\psi^{\text{blip}}_s(c_o) = \psi^{\text{blip}}(c_{o,1} \mathbb{I}(1\in s), \ldots, c_{o,d} \mathbb{I}(d \in s))$, where $c_{o}$ denotes all possibilities of $C_o$. We assume the variation norm of $\psi^{\text{blip}}$ is finite:
\begin{equation*}
\|\psi^{\text{blip}}\|_v = \psi^{\text{blip}}(0) + \sum_{s \subset
 \{1,\ldots,d\}}\int_{0_s}^{\tau_s}|\psi^{\text{blip}}_s (du)| < M.
\end{equation*}
\noindent 
The HAL estimator represents $\psi^{\text{blip}}$ as
\begin{align*}
\psi^{\text{blip}}(c_o) &= \psi^{\text{blip}}(0) + \sum_{s \subset
\{1,\ldots,d\}}\int_{0_s}^{\tau_s}\psi^{\text{blip}}_s(du) \\ &=
\psi^{\text{blip}}(0) + \sum_{s \subset \{1,\ldots,d\}}
\int_{0_s}^{\tau_s} \mathbb{I} (u\leq c_{o,s})\psi^{\text{blip}}_s(du),
\end{align*}
\noindent 
which uses a discrete measure $\psi^{\text{blip}}_m$ with $m$ support
points to approximate this representation. 
For each subset $s$, at time $t=N$, we select as support points the $N$ observed values $\tilde{c}_{o,s}(t)$, $t=1,\ldots,N$, of the context
$C_{o,s}(t)$. 
Then, for each subset $s$, we have a
discrete approximation of $\psi^{\text{blip}}_s$ with support defined by the
actual $N$ observations and point-masses $d_{\psi^{\text{blip}}_{m,s,t}}$, the pointmass
assigned by $\psi^{\text{blip}}_m$ to point $\tilde{c}_{o,s}(t)$, $t=1,\ldots,N$.  
%Note that $\tilde{c}_{o}(t)$ denotes the observed values of $C_o$, $t = 1, \ldots, N$ and for
This approximation consists of a linear combination of basis functions $c_o \rightarrow \phi_{s,t}(c_o)=\mathbb{I}(c_{o,s}\geq
\tilde{c}_{o,s}(t))$ with corresponding coefficients
$d_{\psi^{\text{blip}}_{m,s,t}}$ summed over $s$ and $t=1,\ldots,N$. 

%We define the loss-based dissimilarity for $L_B(\theta)(o,C_o)$ as $\|\psi^{\text{blip}}_n-\psi^{\text{blip}}_0 \|^2_{P_0}$. 
The minimization of the empirical risk $P_n L_B(\theta)(o,C_o)$ of this estimator, $\psi^{\text{blip}}_n$, corresponds to lasso regression with predictors $\phi_{s,t}$ across all subsets $s \subset \{1, \ldots, d\}$ and for $t = 1, \ldots, N$. That is, for
\begin{equation*}
\psi_{\beta}^{\text{blip}}=\beta_0+\sum_{s \subset
    \{1,\ldots,d\}}\sum_{t=1}^N \beta_{s,t}\phi_{s,t}
\end{equation*}
and corresponding subspace $\Psi_{n,M}=\{\psi_\beta:\beta,\beta_0+\sum_{s \subset
\{1,\ldots,d\}}\sum_{t=1}^N |\beta_{s,t}|< M \}$, 
\begin{equation*}
\beta_n=\text{argmin}_{\beta,\beta_0 + \sum_{s \subset \{1,\ldots,d\}}\sum_{t=1}^N| \beta_{s,t}|<M } P_n L_{B_{\beta}}(\theta).
\end{equation*}
    
The linear combination of basis function with non-zero coefficients in the HAL MLE represent a working model. We can use this data adaptively chosen parametric working model  to obtain approximate (non-formal) inference for the blip function. For example, we could use the delta-method to obtain a Wald-type confidence interval for the blip function, recognizing that $\beta_n$ is an MLE for this working model. Alternatively, we use the nonparametric bootstrap, fixing the model that was selected by HAL (to maintain $L1$-norm), and running lasso with the selected model for each bootstrap. We denote the resulting confidence interval bounds around the HAL MLE of the blip as $\pm \text{CI}(\psi_n^{\text{blip}})$, and propose using these bounds around the HAL MLE of the blip - that is, we would let $\pm \text{CI}(\psi_n^{\text{blip}})$ replace $\pm e_t$ in $G_t(x)$. Incorporating the Highly Adaptive Lasso blip estimate into the sampling scheme encourages exploitation of the known uncertainty in the blip estimates so far, allowing for more efficient use of the exploration step then the procedure described in subsection \ref{sotr}.

\section{Targeted Maximum Likelihood Estimator}\label{tmle}

In the following, we build a Targeted Maximum Likelihood Estimator (TMLE) for the target parameter, $\bar{\Psi}(\theta)$ \citep{tmle2006, book2011, book2018}. TML estimation is a multistep procedure, where one first obtains an estimate of the relevant parts of the data-generating distribution using machine learning algorithms and appropriate cross-validation \citep{sl2007, benkeser2018}. The second stage updates the initial fit in a step targeted towards making an optimal bias-variance trade-off for $\bar{\Psi}(\theta)$, instead of the whole density. 

Let $L(\bar{Q})(O(t), C_o(t))$ be a loss function for $\bar{Q}_0$ where $L(\bar{Q}) : \mathcal{O} \times \mathcal{C} \rightarrow \mathbb{R}$; for notational simplicity, we can also write $L(\bar{Q})$, with dependence on $(O(t), C_o(t))$ implied. In particular, we define $L(\bar{Q})$ as the quasi negative log-likelihood loss function, $L(\bar{Q}) = -[Y(t) \log\bar{Q}(C_o(t),A(t)) + (1-Y(t)) \log(1-\bar{Q}(C_o(t),A(t)))]$. The true $\bar{Q}_0$ minimizes the risk under the true conditional density $P_{0,C_0(t)}$:
\begin{equation*}
    P_{0,C_0(t)} L(\bar{Q}_0)(O(t), C_o(t)) = \min_{\bar{Q}} P_{0,C_o(t)} L(\bar{Q})(O(t), C_o(t)).
\end{equation*}
Let $\bar{Q}_N$ be an initial estimator of $\bar{Q}_0$, obtained via online Super Learner and cross-validation suited for dependent data, such as the rolling-window or recursive-origin scheme \citep{Bergmeir2012, benkeser2018}. For a $\bar{Q}_N$ in the statistical model, we define a parameteric working model $\{\bar{Q}_{N,\epsilon} : \epsilon \}$ through $\bar{Q}_N$ with finite-dimensional parameter $\epsilon$; note that  $\bar{Q}_{N,\epsilon = 0} = \bar{Q}_N$. We define a parametric family of fluctuations of the initial estimator $\bar{Q}_N$ of $\bar{Q}_0$ along with the loss function,  $L(\bar{Q})$, so that the linear combination of the components of the derivative of the loss evaluated at $\epsilon=0$ span the efficient influence curve at the initial estimator:

\begin{equation*}
\left \langle \left .  \frac{d}{d\epsilon}L(\bar{Q}_{N,\epsilon}) \right |_{\epsilon =0}\right \rangle \supset D^*_{C_o(t)}(\bar{Q}_N), 
\end{equation*}

\noindent
where we used the notation $\langle S\rangle $ for the linear span of the components of the function $S$. We note that $\{\bar{Q}_{N, \epsilon} : \epsilon \}$ is known as the local least favorable submodel; one could also define a universal least favorable submodel, where the derivative of the loss evaluated at any $\epsilon$ will equal the efficient influence curve at the fluctuated initial estimator $\bar{Q}_{N,\epsilon}$ \citep{onestep2016}. We proceed to maximize the log-likelihood over the parametric model:
\begin{equation*}
    \epsilon_N= \arg\min_{\epsilon} \frac{1}{N} \sum_{t=1}^N
    L(\bar{Q}_{N,\epsilon})(O(t),C_o(t)).
\end{equation*}
In order to perform the update of the conditional expectations, we rely on the logistic fluctuation model,
\begin{equation*}
    \text{logit}(\bar{Q}_{N,\epsilon}) = \text{logit}(\bar{Q}_{N}) + \epsilon H,
\end{equation*}
where $H$ denotes the clever covariate specific to the target parameter, $H = \frac{\mathbb{I}(A(t)=d(C_o(t))}{g_t(A(t)\mid C_o(t))}$.
%We iterate the procedure as $k \rightarrow k+1$, and repeat the process until $\epsilon_N^k \approx 0$. 
The TMLE update, denoted as $\bar{Q}_{N}^*=\bar{Q}_{N,\epsilon_N}$, is the TMLE of $\bar{Q}_0$ which solves the efficient score equation,
\begin{equation*}
    \frac{1}{N}\sum_{t=1}^N D^*(\bar{Q}_N^*)(O(t),C_o(t)) \approx 0.
\end{equation*}
We define the TMLE as the plug-in estimator $\bar{\Psi}(\bar{Q}_N^*)$, obtained by evaluating $\bar{\Psi}$ at the last update of the estimator of $\bar{Q}_0$.

\section{Asymptotic normality of the TMLE}\label{theory}

\subsection{Decomposition of the TMLE estimator}

Our theoretical analysis relies on the fact that the difference between the TML estimator and the target can be decomposed as the sum of (1) the average of a martignale difference sequence, and (2) a martingale process for which we can show an equicontinuity result. We present formally this decomposition in theorem \ref{thm:TMLE_decompositon} below.

\begin{theorem}\label{thm:TMLE_decompositon}
For any $\bar{Q}_1 \in \bar{\mathcal{Q}}$, the difference between the TMLE and its target decomposes as 
\begin{align}
    \bar{\Psi}(\bar{Q}_N^*) - \bar{\Psi}(\bar{Q}_0) = M_{1,N}(\bar{Q}_1) + M_{2,N}(\bar{Q}_N^*, \bar{Q}_1),
\end{align}
with 
\begin{align}
    M_{1,N}(\bar{Q}_1) =& \frac{1}{N} \sum_{t=1}^N D^*(\bar{Q}_1) (C_o(t), O(t)) - P_{0, C_o(t)} D^*(\bar{Q}_1), \\
    M_{2,N}(\bar{Q}^*_N, \bar{Q}_1) =& \frac{1}{N} \sum_{t=1}^N (\delta_{C_o(t), O(t)} - P_{0, C_o(t)}) (D^*(\bar{Q}^*_N) - D^*(\bar{Q}_1)).
\end{align}
\end{theorem}
The first term, ${M}_{1,N}(\bar{Q}_1)$, is the average of a martingale difference sequence, and we will analyze it with a classical martingale central limit theorem. The second term is a martingale process indexed by $\bar{Q} \in \bar{\mathcal{Q}}$, evaluated at $\bar{Q} = \bar{Q}^*_N$. We will prove an equicontinuity result under a complexity condition for a process derived from the function class $\{ D^*(\bar{Q}) : \bar{Q} \in \bar{\mathcal{Q}} \}$, which will imply that if $\bar{Q}^*_N \xrightarrow{P} \bar{Q}_1 \in \bar{\mathcal{Q}}$ then $M_{2,N}(\bar{Q}^*_N, \bar{Q}_1) = o_P(N^{-1/2})$.

\subsection{Analysis of the first term}\label{M1}

A set of sufficient conditions for the asymptotic normality of the term $M_{1,N}(\bar{Q}_1)$ is that (a) the terms $D^*(\bar{Q}_1)(C_o(t), O(t))$ remain bounded, and (b) that the average of the conditional variances of $D^*(\bar{Q}_1)(C_o(t), O(t))$ stabilize. 
A sufficient condition for condition (a) to hold is the following strong version of the positivity assumption.

\begin{assumption}[Strong positivity]\label{assumption:strong_positivity}
There exists $\delta > 0$ such that, for every $t \geq 1$,
\begin{align}
    g_{0,t}(A(t) \mid C_a(t)) \geq \delta, P_0\text{-a.s.}
\end{align}
\end{assumption}

%We state formally below the variance stabilization condition.

\begin{assumption}[Stabilization of the mean of conditional variances]\label{assumption:cond_variances} There exists $\sigma_0^2(\bar{Q}_1) \in (0, \infty)$ such that
\begin{align}
    \frac{1}{N} \sum_{t=1}^N \Var_{Q_0}\left(D^*(\bar{Q}_1)(C_o(t), O(t)) \mid C_o(t) \right) \xrightarrow{d} \sigma_0^2(\bar{Q}_1).
\end{align}
\end{assumption}
\noindent
We formally state below our asymptotic normality result for $M_{1,N}(\bar{Q}_1)$.

\begin{theorem}\label{thm:asymptotic_normality_M1N}
Suppose that assumption \ref{assumption:strong_positivity} and assumption \ref{assumption:cond_variances} hold. Then
\begin{align}
    \sqrt{N} M_{1,N}(\bar{Q}_1) \xrightarrow{d} \mathcal{N}(0, \sigma_0^2(\bar{Q}_1)).
\end{align}
\end{theorem}

\begin{proof}
The result follows directly from various versions of martingale central limit theorems (e.g. theorem 2 in \cite{brown1971}).
\end{proof}

We show in section \ref{section:sufficient_conditions_stab_cond_vars} in the appendix that the conditional variances stabilize under (1) mixing and ergodicity conditions for the sequence $(C_o(t))$ of contexts, and if (2) the design $g_{0,t}$ stabilizes asymptotically. We discuss special cases in which these mixing and ergodicity conditions can be checked explicitly in appendix section \ref{section:sufficient_conditions_stab_cond_vars}. 
We rely on the  empirical variance estimator,
%In section  \ref{section:consistency_empirical_variance_estimator} of the appendix, we give conditions under which the empirical variance estimator,
\begin{align}
    \widehat{\sigma}_N^2 := \frac{1}{N} \sum_{t=1}^N D^*(\bar{Q}_N^*, g_{0,t})^2(C_o(t), O(t)), 
\end{align}
which converges to the asymptotic variance $\sigma_0^2(\bar{Q}_1)$ of $M_{1,N}(\bar{Q}_1).$

\subsection{Negligibility of the second term}

In this susbsection, we give a brief overview of the analysis of the term $M_{2,N}(\bar{Q}_N^*, \bar{Q}_1)$, which we carry out in detail in appendix section \ref{section:analysis_martingale_process_term}.
We show that $M_{2,N}(\bar{Q}_N^*, \bar{Q}_1) = o_P(N^{-1/2})$ by proving  an equicontinuity result for the process $\{M_{2,N}(\bar{Q}, \bar{Q}_1) : \bar{Q} \in \bar{\mathcal{Q}} \}$. Our equicontinuity result relies on a measure of complexity for the process 
\begin{align}
    \Xi_N := \left\lbrace \left(D^*(\bar{Q}, g_{0,t})(C_o(t), O(t)) - D^*(\bar{Q}_1, g_{0,t})(C_o(t), O(t)) \right)_{t=1}^N : \bar{Q} \in \bar{\mathcal{Q}} \right\rbrace, \label{eq:def_can_gdt_process}
\end{align}
which we refer to as \textit{sequential bracketing entropy}, introduced by \cite{vandeGeer2000} for the analysis of martingale processes. We relegate the formal definition of the \textit{sequential bracketing entropy} to the appendix section \ref{section:analysis_martingale_process_term}. In particular, we denote $N_{[\,]}(\epsilon, b, \Xi_N, \bar{O}(N))$ as the \textit{sequential bracketing number} of $\Xi_N$ corresponding to brackets of size $\epsilon$. Our equicontinuity result is a sequential equivalent of similar results for i.i.d. settings (e.g. \cite{vanderVaart&Wellner96}) and similarly relies on a Donsker-like condition.
%(we refer to section \ref{section:analysis_martingale_process_term} for the meaning $b$ and $\bar{O}(N)$). 

\begin{assumption}[Sequential Donsker condition]\label{assumption:Donsker_condition_EIF_process}
Define the sequential bracketing entropy integral as $J_{[\,]}(\epsilon, b, \Xi_N, \bar{O}(N)) := \int_0^\epsilon \sqrt{\log (1 + \mathcal{N}_{[\,]}(u, b, \Xi_N, \bar{O}(N))} du.$
Suppose that there exists a function $a : \mathbb{R}^+ \to \mathbb{R}^+$ that converges to $0$ as $\delta \to 0$, such that 
\begin{align}
    J_{[\,]}(\epsilon, b, \Xi_N, \bar{O}(N)) \leq a(\delta).
\end{align}
\end{assumption}
\noindent
Note that a sufficient condition for assumption \ref{assumption:Donsker_condition_EIF_process} to hold is that $\log (1 + \mathcal{N}_{[\,]}(u, b, \Xi_N, \bar{O}(N)) \leq C \epsilon^{-p}$, with $p \in (0,2)$ and $C > 0$ a constant that does not depend on $N$.

\begin{assumption}[$L_2$ convergence of the outcome model]\label{assumption:convergence_outcome_model} It holds that
$\| \bar{Q}^*_N - \bar{Q}_1 \|_{2, g^*, h_N} = o_P(1)$, where $h_N$ is the empirical measure $h_N:= N^{-1} \sum_{t=1}^N \delta_{C_o(t)}$.
\end{assumption}

\begin{theorem}[Equicontinuity of the martingale process term]\label{thm:equicontinuity_M2N}
Consider the process $\Xi_N$ defined in equation \eqref{eq:def_can_gdt_process}. Suppose that assumptions \ref{assumption:strong_positivity}, \ref{assumption:Donsker_condition_EIF_process} and \ref{assumption:convergence_outcome_model} hold. Then $M_{2,N}(\bar{Q}^*_N, \bar{Q}_1) = o_P(N^{-1/2})$.
\end{theorem}

\subsection{Asymptotic normality theorem}

As an immediate corollary of theorems \ref{thm:asymptotic_normality_M1N} and \ref{thm:equicontinuity_M2N}, we have the following asymptotic normality result for our TML estimator.

\begin{theorem}[Asymptotic normality of the TMLE]\label{thm:asymptotic_normality_TMLE}
Suppose that assumptions \ref{assumption:strong_positivity}, \ref{assumption:cond_variances}, \ref{assumption:Donsker_condition_EIF_process} and \ref{assumption:convergence_outcome_model} hold. Then
\begin{align}
    \sqrt{N} \left( \bar{\Psi}(\bar{Q}_N^*) - \bar{\Psi}(\bar{Q}_0) \right) \xrightarrow{d} \mathcal{N}(0, \sigma_0^2(\bar{Q}_1)).
\end{align}
\end{theorem}

\noindent
%As mentioned earlier, under conditions we discuss in appendix section \ref{section:consistency_empirical_variance_estimator}, 
The empirical variance estimator $\widehat{\sigma}_N^2$ converges in probability to $\sigma_0^2(\bar{Q}_1)$, which implies that 
\begin{align}
    \widehat{\sigma}_N^{-1} \sqrt{N}  \left( \bar{\Psi}(\bar{Q}_N^*) - \bar{\Psi}(\bar{Q}_0) \right) \xrightarrow{d} \mathcal{N}(0, 1).
\end{align}
Therefore, denoting $q_{1-\alpha/2}$ the $1-\alpha/2$-quantile of the standard normal distribution, we have that the confidence interval
\begin{align}
    \left[\bar{\Psi}(\bar{Q}_N^*) - \frac{q_{1-\alpha/2} \widehat{\sigma}_N}{\sqrt{N}},  \bar{\Psi}(\bar{Q}_N^*) + \frac{q_{1-\alpha/2} \widehat{\sigma}_N}{\sqrt{N}}\right]
\end{align}
has asymptotic coverage $1-\alpha$ for the target $\bar{\Psi}(\bar{Q}_0)$.

\section{Simulations}\label{simulations}

\noindent
In this section we present simulation results concerning the adaptive learning of the optimal individualized treatment rule estimated using machine learning methods for a single time-series. We focus on the stochastic sampling scheme described in subsection \ref{sotr}, and explore performance of our estimator with different initial sample sizes and consequent sequential updates. We consider binary outcome and treatment, but note that the results will be comparable for continuous bounded outcome. Finally, unless specified otherwise, we present coverage of the mean under the current estimate of the optimal individualized treatment rule at each update based on 500 Monte Carlo draws. We set the reference treatment mechanism to a balanced design, assigning treatment with probability $0.5$ for the data draw used to learn the initial estimate of the optimal individualized treatment rule.

\subsection{Simulation 1a}\label{sim1a}

We explore a simple dependence setting first, emphasising the connection with i.i.d sequential settings. We data consists of a binary treatment ($A(t) \in \{0,1\}$) and outcome ($Y(t) \in \{0,1\}$). The time-varying covariate $W(t)$ decomposes as $W(t) \equiv (W_1(t), W_2(t))$ with binary $W_1$ and continuous $W_2$. The outcome $Y$ at time $t$ is conditionally drawn given $\{A(t), Y(t-1), W_1(t-1)\}$ from a Bernoulli distribution, with success probability defined as $1.5*A(t) + 0.5 * Y(i-1) - 1.1*W_1(i-1)$. We generate the initial sample of size $t=1000$ and $t=500$ by first drawing a set of four $O(t)$ samples randomly from binomial and normal distributions in order to have a starting point to initiate time dependence. After the first 4 draws, we draw $A(t)$ from a binomial distribution with success probability 0.5, $Y(t)$ from a Bernoulli distribution with success probability dependent on $\{A(t), A(t-1), Y(t-1), W_2(t-1)\}$, followed by $W_1(t)$ conditional on $\{Y(t-1),W_1(t-1),W_2(t-1)\}$ and $W_2(t)$ conditional on $\{A(t-1), Y(t-1), W_1(t-1)\}$. After $t=1000$ or $t=500$, we continue to draw $O(t)$ as above, but with $A(t)$ drawn from a stochastic intervention approximating the current estimate $d_{\bar{Q}_{t-1}}$ of the optimal rule $d_{\bar{Q}_0}$. This procedure is repeated until reaching a specified final time point indicating the end of a trial. Our estimator of $\bar{Q}_{0}$, and thereby the optimal rule $d_0$, is based on an online super-learner with an ensemble consisting of multiple algorithms, including simple generalized linear models, penalized regressions, HAL and extreme gradient boosting \citep{coyle2018sl3}.  For cross-validation, we relied on the online cross-validation scheme, also known as the recursive scheme in the time-series literature. The sequences $\{c_t\}_{t \geq 1}$ and $\{e_t\}_{t \geq 1}$ are chosen constant, with $c_{\infty} = 10\%$ and $e_{\infty} = 5\%$. The TMLEs are computed at sample sizes a multiple of 200, and no more than 1800 (for initial $t=1000$) or 1300 (for initial $t=500$), at which point sampling is stopped. We use the coverage of asymptotic 95$\%$ confidence intervals to evaluate the performance of the TMLE in estimating the average across time $t$ of the $d_{\bar{Q}_{t-1}}$-specific mean outcome. The exact data-generating distribution used is as follows:
\noindent
\begin{small}
\begin{align*}
A(0:4) &\sim \text{Bern}(0.5) \\
Y(0:4) &\sim \text{Bern}(0.5) \\
W_1 (0:4) &\sim \text{Bern}(0.5) \\
W_2 (0:4) &\sim \text{Normal}(0,1) \\
A(4:t) &\sim \text{Bern}(0.5) \\
Y(4:t) &\sim \text{Bern}(expit(1.5*A(i) + 0.5*Y(i-1) - 1.1*W_1(i-1) )) \\
W_1 (4:t) &\sim \text{Bern}(expit(0.5*W_1(i-1) - 0.5*Y(i-1) + 0.1*W_2(i-1))) \\
W_2 (4:t) &\sim \text{Normal}(0.6*A(i-1) + Y(i-1) - W_1(i-1), sd=1) \\
A(t:1800) &\sim d_{\bar{Q}_{t-1}} \\
Y(t:1800) &\sim \text{Bern}(expit(1.5*A(i) + 0.5*Y(i-1) - 1.1*W_1(i-1) )) \\
W_1 (t:1800) &\sim \text{Bern}(expit(0.5*W_1(i-1) - 0.5*Y(i-1) + 0.1*W_2(i-1))) \\
W_2 (t:1800) &\sim \text{Normal}(0.6*A(i-1) + Y(i-1) - W_1(i-1), sd=1).
\end{align*}
\end{small}
\noindent
From Table \ref{table:table1}, we can see that the 95$\%$ coverage for the average across time of the counterfactual mean outcome under the current estimate of the optimal dynamic treatment approaches nominal coverage with increasing time-steps, for both $t=500$ and $t=1000$ length of the initial time-series. The mean conditional variance stabilizes with increasing time-steps, as illustrated in Table \ref{table:table2} and Figure \ref{fig1}A, thus satisflying assumption \ref{assumption:cond_variances} necessary for showing asymptotic normality of the TML estimator. 

\subsection{Simulation 1b}\label{sim1b}

\begin{figure}[p]
\centering
     \includegraphics[scale=0.75]{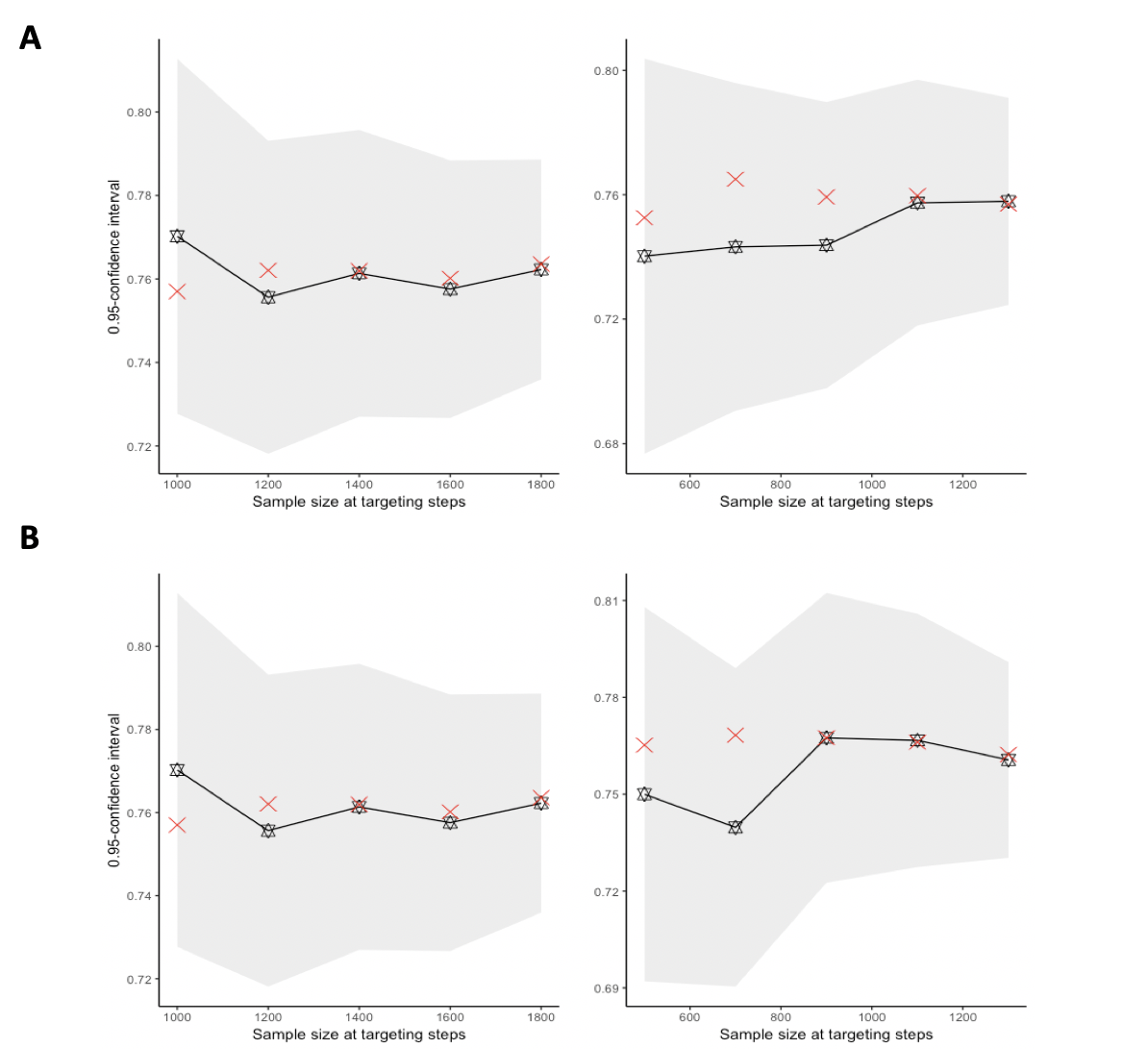}
    \caption{Illustration of the data-adaptive inference of the mean reward under the optimal treatment rule with initial sample size $n=1000$ and $n=500$ for Simulation 1a and 1b. The red crosses reflect successive values of the data-adaptive true parameter, with stars representing the estimated parameter with the corresponding $95\%$ confidence interval for the data-adaptive parameter.}
    \label{fig1}
\end{figure}

In Simulation 1b, we explore the behavior of our estimator in case of more elaborate dependence. As in Simulation 1a, we only consider binary treatment ($A(t) \in \{0,1\}$) and outcome ($Y(t) \in \{0,1\}$), with binary and continuous time-varying covariates. We set the reference treatment mechanism to a balanced treatment mechanism assigning treatment with probability $P(A(t) = 1) = 0.5$, and generate the initial sample of size $t=(1000, 500)$ by sequentially drawing $W_1(t), W_2(t), A(t), Y(t)$. As before, upon the first $t=1000$ or $t=500$ time-points, we continue to draw $O(t)$ with $A(t)$ sampled from a stochastic intervention approximating the current estimate $d_{\bar{Q}_{t-1}}$ of the optimal rule $d_{\bar{Q}_0}$. The estimator of the optimal rule  $d_{\bar{Q}_0}$ was based on an ensemble of machine learning algorithms and regression-based algorithms, with honest risk estimate achieved by utilizing online cross-validation scheme with validation set size of 30. The sequences $\{c_t\}_{t \geq 1}$ and $\{e_t\}_{t \geq 1}$ were set to $10\%$ and $ 5\%$, respectively. The TMLEs are computed at initial $t=1000$ or $t=500$, and consequently at sample sizes being a multiple of 200, and no more than 1800 (or 1300), at which point sampling is stopped. The exact data-generating distribution used is as follows: 
\noindent
\begin{small}
\begin{align*}
A(0:4),Y(0:4), W_1 (0:4) &\sim \text{Bern}(0.5) \\
%Y(0:4) &\sim \text{Bern}(0.5) \\
%W_1 (0:4) &\sim \text{Bern}(0.5) \\
W_2 (0:4) &\sim \text{Normal}(0,1) \\
A(4:t) &\sim \text{Bern}(0.5) \\
Y(4:t) &\sim \text{Bern}(expit(1.5*A(i) + 0.5*Y(i-3) - 1.1*W_1(i-4) )) \\
W_1 (4:t) &\sim \text{Bern}(expit(0.5*W_1(i-1) - 0.5*Y(i-1) + 0.1*W_2(i-2))) \\
W_2 (4:t) &\sim \text{Normal}(0.6*A(i-1) + Y(i-1) - W_1(i-2), sd=1) \\
A(t:1800) &\sim d_{\bar{Q}_{t-1}} \\
Y(t:1800) &\sim \text{Bern}(expit(1.5*A(i) + 0.5*Y(i-3) - 1.1*W_1(i-4) )) \\
W_1 (t:1800) &\sim \text{Bern}(expit(0.5*W_1(i-1) - 0.5*Y(i-1) + 0.1*W_2(i-2))) \\
W_2 (t:1800) &\sim \text{Normal}(0.6*A(i-1) + Y(i-1) - W_1(i-2), sd=1).
\end{align*}
\end{small}
\noindent
As demonstrated in Table \ref{table:table1}, the TML estimator approaches 95$\%$ coverage with increasing number of time points with more elaborate dependence structure as well. The assumption of stabilization of the mean of conditional variances is shown to be valid in Table \ref{table:table2} and Figure \ref{fig1}B, allowing for the asymptotic coverage $1-\alpha$ for the target $\bar{\Psi}(\bar{Q}_0)$.

\begin{table}[p]
\centering
\begin{tabular}{@{}lrrrrrr@{}}\toprule 
& \textbf{$t$} & $\textbf{Cov}_{t}$ & $\textbf{Cov}_{t_1}$ & $\textbf{Cov}_{t_2}$ & $\textbf{Cov}_{t_3}$  & $\textbf{Cov}_{t_4}$\\ \midrule

\textbf{Simulation 1a} & 1000  & 92.60 & 94.00 & 95.20 & 95.40 & 95.80\\  
\textbf{Simulation 1a} & 500  & 90.00 & 93.20 & 93.80 & 94.80 & 94.60\\  \hline

\textbf{Simulation 1b} & 1000  & 92.60 & 92.60 & 93.00 & 93.40 & 93.80 \\ 
\textbf{Simulation 1b} & 500  & 89.60 & 90.20 & 89.90 & 90.80 & 91.40\\  \hline
\bottomrule
\end{tabular}

\caption{The 95$\%$ coverage for the average across time of the counterfactual mean outcome under the current estimate of the optimal dynamic treatment at time points $t$, $t_1 = t+200$, $t_2 = t+400$, $t_3 = t+600$ and $t_4 = t+800$. The first $t$ time points sample treatment with probability 0.5. The sequences $\{c_n\}_{t \geq 1}$ and $\{e_n\}_{t \geq 1}$ are chosen constant, with $c_{\infty} = 10\%$ and $e_{\infty} = 5\%$. TMLEs are computed at $t = \{500, 1000\}$, $t_1$, $t_2$, $t_3$ and $t_4$, with sequential updates being of size 200. The results are reported over 500 Monte-Carlo draws for Simulations 1a and 1b with initial sample sizes 1000 and 500.}
\label{table:table1}
\end{table}

\vspace{5mm}
\begin{table}[p]  \centering
\begin{tabular}{@{}lrrrrrr@{}}\toprule
& \textbf{$t$} & $\textbf{Var}_{t}$ & $\textbf{Var}_{t_1}$ & $\textbf{Var}_{t_2}$ & $\textbf{Var}_{t_3}$  & $\textbf{Var}_{t_4}$\\ \midrule

\textbf{Simulation 1a} & 1000  & 0.0018 & 0.0019 & 0.0017 & 0.0016 & 0.0004 \\  
\textbf{Simulation 1a} & 500  & 0.0011 & 0.0024 & 0.0035 & 0.0014 & 0.0011 \\  \hline

\textbf{Simulation 1b} & 1000  & 0.0072 & 0.0075 & 0.0069 & 0.0067 & 0.0018 \\  
\textbf{Simulation 1b} & 500  & 0.0199 & 0.0171 & 0.0187 & 0.0152 & 0.0087 \\  \hline
\bottomrule
\end{tabular}

\caption{Variance for the average across time of the counterfactual mean outcome under the current estimate of the optimal dynamic treatment at time points $t$, $t_1 = t+200$, $t_2 = t+400$, $t_3 = t+600$ and $t_4 = t+800$, over 500 Monte-Carlo draws for Simulations 1a and 1b with initial sample sizes 1000 and 500.}
\label{table:table2}
\end{table}

\section{Conclusions}\label{conclusions}

In this manuscript, we consider causal parameters based on observing a single time series with asymptotic results derived over time $t$. The data setup constitutes a typical longitudinal data structure, where within each $t$-specific time-block one observes treatment and outcome nodes, and possibly time-dependent covariates in-between treatment nodes. Each $t$-specific data record $O(t)$ is viewed as its own experiment in the context of the observed history $C_o(t)$, carrying information about a causal effect of the treatment nodes on the next outcome node. While in this work we concentrate on single time point interventions, we emphasize that our setup can be easily generalized to context specific causal effects of multiple time point interventions, therefore estimating the causal effect of $A(t:t+k)$ on future $Y(t+k)$. 

A key assumption necessary in order to obtain the presented results is that the relevant history for generating $O(t)$, given the past $\bar{O}(t-1)$, can be summarized by a fixed dimensional summary $C_o(t)$. We note that our conditions allow for $C_o(t)$ to be a function of the whole observed past, allowing us to avoid Markov-order type assumptions that limit dependence on recent, or specifically predefined past. Components of $C_o(t)$ that depend on the whole past, such as an estimate of the optimal treatment rule based on $(O(1),\ldots,O(t-1))$, will typically converge to a fixed function of a recent past - so that the martingale condition on the stabilization of the mean of conditional variances holds.
%$\frac{1}{N} \sum_t P_{\theta_0,C_o(t)}\{D^*_{C_o(t)}(\theta^*)\}^2\rightarrow \sigma^2_0$ will still hold. 

Due to the dimension reduction assumption, each $t$-specific experiment corresponds to drawing from a conditional distribution of $O(t)$ given $C_o(t)$. We assume that this conditional distribution is either constant in time or is parametrized by a constant function. As such, we can learn the true mechanism that generates the time-series, even when the model for the mechanism is nonparametric. With the exception of parametric models allowing for maximum likelihood estimation, we emphasize that statistical inference for proposed target parameters of the time-series data generating mechanism is a challenging problem which requires targeted machine learning. 

The work of \cite{vanderLaan2018onlinets} and \cite{kallus2019efficiently} studies marginal causal parameters, marginalizing over the distribution of $C_o(t)$, defined on the same statistical model as the parameter we consider in this article. In particular, \cite{vanderLaan2018onlinets} define target parameters and estimation of the counterfactual mean of a future (e.g., long term) outcome under a stochastic intervention on a subset of the treatment nodes, allowing for extensions to single unit causal effects. As such, the target parameter proposed by \cite{vanderLaan2018onlinets} addresses the important question regarding the distribution of the outcome at time $t$, had we intervened on some of the past treatment nodes in a (possibly single) time-series. While important, the TMLE of such target parameters are challenging to implement due to their reliance on the density estimation of the marginal density of $C_o(t)$ averaged across time $t$. Additionally, we remark that such marginal causal parameters cannot be robustly estimated if treatment is sequentially randomized, due to the lack of double robustness of the second order remainder. 

In this work, we focus on a context-specific target parameter is order to explore robust statistical inference for causal questions based on observing a single time series of a particular unit. We note that for each given $C_o(t)$, any intervention-specific mean outcome $EY_{g^*}(t)$ with $g^*$ being a stochastic intervention w.r.t. the conditional distribution of $P_{C_o(t)}$ (with deterministic rule being a special case), represents a well studied statistical estimation problem based on observing many i.i.d. copies. Even though we do not have repeated observations from the $C_o(t)$-specific distribution at time $t$, the collection $(C_o(t),O(t))$ across all time points represent the analogue of an i.i.d. data set $(C_o(t),O(t))\sim_{iid} P_0$, where $C_o(t)$ can be viewed as a baseline covariate for the longitudinal causal inference data structure; we make the connection with the i.i.d. sequential design in one of our simulations. The initial estimation step of the TMLE should still respect the known dependence in construction of the initial estimator, by relying on appropriate estimation techniques developed for dependent data. Similarly, variance estimation can proceed as in the i.i.d case using the relevant i.i.d. efficient influence curve. This insight relies on the fact that the TMLE in this case allows for the same linear approximation as the TMLE for i.i.d. data, with the martingale central limit theorem applied to the linear approximation instead. Since the linear expansion of the time-series TMLE for context-specific parameter is an element of the tangent space of the statistical model, our derived TMLE is asymptotically efficient.

%To emphasize the importance of our work in applied settings, we provide an exciting application of the context-specific parameter in the settings where the optimal individualized rule is learned adaptively from a single observed time-series. 

Our motivation for studying the proposed context-specific parameter strives from its important role in precision medicine, in which one wants to tailor the treatment rule to the individual observed over time. In particular, we derive a TMLE which uses only the past data $\bar{O}(t-1)$ of a single unit in order to learn the optimal treatment rule for assigning $A(t)$ to maximize the mean outcome $Y(t)$. Here, we assign the treatment at the next time point $t+1$  according to the current estimate of the optimal rule, allowing for the time-series to learn and apply the optimal treatment rule at the same time. The time-series generated by the described adaptive design within a single unit can be used to estimate, and most importantly provide inference, for the average across all time-points $t$ of the counterfactual mean outcome of $Y(t)$ under the estimate $d(C_o(t))$ of the optimal rule at a relevant time point $t$. Assuming that the estimate of the optimal rule is consistent, as the number of time-points increases, our target parameter converges to the mean outcome one would have obtained had they carried out the optimal rule from the start. As such, we can effectively learn the optimal rule and simultaneously obtain valid inference for its performance. Interestingly, this does not provide inference relative to, for example, the control that always assigns $A(t)=0$. This is due to the fact that by assigning treatment $A(t)$ according to a rule, the positivity assumption needed to learn $\frac{1}{N}\sum_t E(Y_{A(t)=0}(t)\mid C_o(t))$ is violated. However, we note that one can safely conclude that one will not be worse than this control rule, even when the control rule is equal to the optimal rule. If one is interested in inference for a contrast based on a single time-series, then we advocate for random assignment between the control and estimate of optimal rule. As such, our proposed methodology still allows to learn the desired contrast. 

Finally, we note that while the context-specific parameter enjoys many important statistical and computational advantages as opposed to the marginal target parameter based on a single time-series, the formulation employed in this article is only sensible if one is interested in the causal effect of treatment on a short-term outcome. In particular, if the amount of time necessary to collect outcome $Y(t)$ in $O(t)$ is long, then generating a long time series would take too much time to be practically useful. If one is interested in causal effects on a long term outcome and is willing to forgo utilizing known randomization probabilities for treatment, we advocate for the marginal target parameters as described in previous work by \cite{vanderLaan2018onlinets} or \cite{kallus2019efficiently}. 

\section{Acknowledgments}
Research reported in this publication was supported by the National Institute Of Allergy And Infectious Diseases of the National Institutes of Health under Award Number R01AI074345. The content is solely the responsibility of the authors and does not necessarily represent the official views of the National Institutes of Health.

\newpage
\section{Appendix}

\noindent

\subsection{Comparison with marginal parameters}\label{subsection:comparison_marg_params}

We present below two alternative statistical parameters defined on the same statistical model as the parameter we consider in this article, and which were considered in previous works \citep{vanderLaan2018onlinets, kallus2019efficiently}. The parameters are \textit{marginal}, as opposed to context-specific parameters we consider in the present article. The definition of the marginal parameters entails integrating against certain marginal distributions of contexts, as we make explicit below.

%\paragraph{G-computation formula.} 
Consider the distribution $P_{Q,g^*}$ over infinite sequences taking values in the infinite cartesian product space $\times_{t=1}^\infty \mathcal{O}$, defined from the factors of $P \in \mathcal{M}$ by the following G-computation formula:
\begin{align}
    P_{Q,g^*}\left((o(t))_{t=1}^\infty\right) := P_{C_o(1)}(c_o(1)) \prod_{t=1}^\infty g^*(a(t) \mid c_o(t)) Q(y(t) \mid c_o(t)) Q_w(w(t) \mid c_o(t)).
\end{align}
Let $(O^*(t))_{t=1}^\infty \sim P_{Q,g^*}$, with $O^*(t) = (A^*(t), Y^*(t), W^*(t))$.

\subsubsection{Marginal parameter by van der Laan et al. (2018)}

As a first example of a marginal parameter, \cite{vanderLaan2018onlinets} consider a class of parameters which includes
\begin{align}
    \Psi_{1,\tau}(P) := E_{Q,g^*}[Y^*(\tau)],
\end{align}
for $\tau \geq 1$. Under the causal identifiability assumptions 
\ref{assumption:sequential_randomization} and \ref{assumption:positivity}, $\Psi_{1,\tau}(P_0)$ equals the mean outcome we would obtain at time $\tau$, under a counterfactual time series with initial context distribution $P_{0,C_o(1)}$ and intervention $g^*$ (instead of the observed intervention $g$) at every time point. We note that $P_{0,C_o(1)}$ is the initial, observed data-generating distribution. The canonical gradient of $\Psi_{1,\tau}$ w.r.t. our model $\mathcal{M}$ (where $\mathcal{M}$ assumes $P_{C_o(1)}$ known\footnote{If we instead supposed that $P_{C_o(1)}$ is unknown and lies in a certain model $\mathcal{M}_{P_{C_o(1)}}$, the canonical gradient would have one additional component, which would be lying in the tangent space of $\mathcal{M}_{P_{C_o(1)}}$. As far as the conditional parameter of the main text are concerned, this distinction has no effect, as these do not depend on the marginal distribution of contexts and therefore its canonical gradient has no components in the tangent spaces corresponding to the context distributions.}) is
\begin{align}
    D^*(P)(o^N) := \frac{1}{N} \sum_{t=1}^N \bar{D}(Q,\omega,g)(c_o(t), o(t))
\end{align}
with
\begin{align}
    \bar{D}(Q,\omega,g)(c_o,o) := \sum_{s=1}^\tau \omega_s(c) \frac{g^*(a\mid c_o)}{g(a \mid c_o)}& \left\lbrace  E_{Q,g^*}[Y^*(\tau) \mid O^*(s) = o, C_o^*(s) = c_o] \right. \\
    & \left.- E_{Q,g^*} [ Y^*(\tau) \mid A^*(s) =a, C_o^*(s) =c_o] \right\rbrace,
\end{align}
with $\omega_s(c_o) = h_{C^*_o(s)}(c_o) / \bar{h}_N(c_o)$, where
\begin{align}
    h_{C_o(s)}(c_o) =& P_{Q,g}[C_o(s) = c_o], \\
    \bar{h}_{N}(c_o) =& \frac{1}{N} \sum_{t=1}^N h_{C_o(t)}(c_o), \\
    \text{and } h_{C^*_o(s)}(c_o) =& P_{Q,g^*}[C^*_o(s) = c_o]
\end{align}
are the marginal density of context $C_o(s)$ under $P$, the average thereof over observed time points $t = 1,\ldots, N$, and the marginal density of context $C^*_o(s)$ under $P_{Q,g^*}$. We note that $\Psi_{1,1}$ is the marginal equivalent of our parameter $\Psi_{C_o(1)}$. Specifically, $\Psi_{1,1}(P) = \int dP_{C_o(1)}(c_o(1)) \Psi_{c_o(1)}(P)$.

\subsubsection{Marginal parameter by Kallus and Uehara (2019)}

Let $\gamma \in (0,1)$. \cite{kallus2019efficiently} consider the parameter
\begin{align}
    \Psi_2(P) := & E_{Q,g^*} \left[ \sum_{\tau=1}^\infty \gamma^\tau Y^*(\tau) \right] \\
    =& \sum_{\tau \geq 1} \gamma^\tau \Psi_{1,\tau}(P).
\end{align}

Under the causal identifiability assumptions \ref{assumption:sequential_randomization} and \ref{assumption:positivity}, $\Psi_2(P_0)$ is the expected total discounted outcome from time point $1$ until $\infty$ that we would get if we carried out intervention $g^*$ forever - starting from initial context distribution $P_{0,C_o(1)}$ as in the observed data generating distribution. The canonical gradient $\Psi_2$ w.r.t. $\mathcal{M}$ (again, supposing that $\mathcal{M}$ considers $P_{0,C_o(1)}$ known) is
\begin{align}
    D^*(P)(o^N) := \frac{1}{N} \sum_{t=1}^N \bar{D}(Q, \omega, g)(c_o(t), o(t)),
\end{align}
with
\begin{align}
    \bar{D}(Q,\omega, g)(c_o, o) := \sum_{s=1}^\infty \omega_s(c_o) \frac{g^*(a \mid c_o)}{g(a \mid c_o)} \left\lbrace y + \gamma V_{1,Q,g^*}(c_o, o) - V_{2,Q,g^*}(c_o, a)\right\rbrace,
\end{align}
where $\omega_s$ is defined as in the previous example, and 
\begin{align}
    V_{1,Q,g^*}(c_o,o) :=& E_{Q,g^*} \left[ \sum_{\tau \geq 2} \gamma^\tau Y^*(\tau) \mid C^*_o(1) = c_o, O^*(1) =o \right] \\
    \text{and } V_{2,Q,g^*}(c_o,o) :=& E_{Q,g^*} \left[ \sum_{\tau \geq 1} \gamma^\tau Y^*(\tau) \mid A^*(1) = a, O^*(1) =o \right].
\end{align}

\subsubsection{Robustness properties}

In this article we are concerned with adaptive trials where the intervention is controlled by the experimenter, hence $g_0$ is known; we therefore only consider the case $g = g_0$. Under $g = g_0$, both parameters $\Psi' \in \{\Psi_{1,\tau}, \Psi_2 \}$ defined above admit a first order expansion of the form
\begin{align}
    \Psi'(P) - \Psi'(P_0) = - P_0 D^*(P) + R'(Q,Q_0, \omega, \omega_0),
\end{align}
where $R'$ is a second-order remainder term such that $R(Q,Q_0, \omega, \omega_0)=0$ if either $Q = Q_0$ and $\omega = \omega_0$. While this resembles a traditional double-robustness property, as that which holds in the i.i.d. setting for the ATE or in the time series setting for our conditional parameter (as opposed to arbitrary time-series dependence or Markov decision process) it is important to note the following:
\begin{enumerate}
    \item For $\Psi' \in \{\Psi_{1,\tau}, \Psi_2 \}$, knowledge of the treatment mechanism is not sufficient to guarantee that the remainder term is zero; we direct the interested reader to \cite{vanderLaan2018onlinets} for the exact form of $R'$.
    \item The parameters $\omega$ and $Q$ are not variation independent, as appears explicitly from the definition of $\omega_s$. In fact, when estimating $\omega_s$ from a single time series, one must a priori rely on an estimator of $Q$ to obtain estimates of $\omega_s$ (see \cite{vanderLaan2018onlinets}). Therefore, if the estimator of $Q$ is inconsistent, the corresponding estimator of $\omega_s$ will be inconsistent as well.
\end{enumerate}

\subsection{Sufficient conditions for the stabilization of conditional variances}\label{section:sufficient_conditions_stab_cond_vars}

Assumption \ref{assumption:cond_variances} on the stabilization of the conditional variance of the canonical gradient can be checked under mixing conditions on the sequence of context $(C_o(t))$, and under the condition that the design $g_{0,t}$ converges to a fixed design. We state formally below such a set of conditions.

\begin{assumption}[Convergence of the marginal law of contexts]\label{assumption:context_marginal_convergence} Suppose that the marginal law of contexts converges to a limit law, that is $C_o(t) \xrightarrow{d} C_\infty$, for some random variable $C_\infty$.
\end{assumption}

The next assumption is a mixing condition in terms of $\rho$-mixing. We first recall the notion of $\rho$-mixing.

\begin{definition}[$\rho$-mixing]\label{def:rho-mixing} Consider a couple of random variables $(Z_1, Z_2) \sim P$. The $\rho$-mixing coefficient, or maximum correlation coefficient of $Z_1$ and $Z_2$ is defined as
\begin{align}
    \rho_{P}(Z_1, Z_2) := \sup \left\lbrace \Corr(f_1(Z_1), f_2(Z_2)) : f_1 \in L_2(P_{Z_1}), f_2 \in L_2(P_{Z_2})\right\rbrace.
\end{align}

\end{definition}

\begin{assumption}[$\rho$-mixing condition]\label{assumption:rho_mixing} Suppose that
\begin{align}
    \sup_{t \geq 1} \sum_{s =1}^N  \rho(C_o(t), C_o(t+s)) = o(N)
\end{align}
\end{assumption}

Observe that if $g$ is common across time points, the process $(C_o(t))$ is an homogeneous Markov chain. Conditions under which homogeneous Markov chains have marginal law converging to a fixed law and are mixing have been extensively studied. A textbook example, albeit perhaps a bit too contrived for many specifications of the setting of our current article, is when the Markov chain has finite state space and the probability of transitioning between any two states from one time point to the next is non-zero. In this case, ergodic theory shows that the transition kernel of the Markov chain admits a so-called invariant law - the marginal laws converge exponentially fast (in total variation distance) to the invariant law, and the mixing coefficients have finite sum. We refer the interested reader to the survey paper by \cite{bradley2005} for more general conditions under which Markov chains have convergent marginal laws and are strongly mixing (for various types of mixing coefficients, one of them being $\rho$-mixing)

\begin{assumption}[Design stabilization]\label{assumption:design_stabilization}
    There is a design $g_\infty$ such that $\|g_{0,t} - g_\infty \|_{1, P_{g^*, h_{0,t}}} = o(1)$, and $g_\infty \geq \delta$, for some $\delta > 0.$
\end{assumption}
\noindent
We note that, as we will always use assumption \ref{assumption:design_stabilization} along with assumption \ref{assumption:strong_positivity}, we will suppose that the constant $\delta$ in the statement of both assumptions is the same.

\begin{lemma}[Conditional variance stabilization under mixing]\label{lemma:cond_var_stab_under_mixing} Suppose that assumptions \ref{assumption:strong_positivity}, \ref{assumption:context_marginal_convergence} and \ref{assumption:rho_mixing} hold. Then assumption \ref{assumption:cond_variances} holds.
\end{lemma}
\noindent
We dedicate the appendix subsection \ref{lemma:cond_var_stab_under_mixing_proof} to the proof of lemma \ref{lemma:cond_var_stab_under_mixing}. 

%\paragraph{Discussion of assumptions \ref{assumption:context_marginal_convergence} and \ref{assumption:rho_mixing}.} 

%\subsection{Consistency of the asymptotic variance estimator}\label{section:consistency_empirical_variance_estimator}

\subsection{Analysis of the martingale process term}\label{section:analysis_martingale_process_term}

We analyze the martingale process $\{M_{2,N}(\bar{Q}, \bar{Q}_1) : \bar{Q} \in \bar{\mathcal{Q}} \}$ under a measure of complexity introduced by \cite{vandeGeer2000}, which we will refer to in the present work as \textit{sequential bracketing entropy}. We state below the definition of sequential bracketing entropy particularized to our setting.

\begin{definition}[Sequential bracketing entropy]\label{def:seq_bracket_ent} Consider a stochastic process of the form $\Xi_N := \{ (\xi_t(f))_{t=1}^N : f \in \mathcal{F} \}$ where $\mathcal{F}$ is an index set such that, for every $f \in \mathcal{F}$, $t \in [N]$, $\xi_t(f)$ is an $\bar{O}(t)$-measurable real valued random variable. We say that a collection of random variables of the form $\mathcal{B} := \{ (\Lambda_t^j, \Upsilon_t^j)_{t=1}^N : j \in [J] \}$ is an $(\epsilon, b, \bar{O}(N))$ bracketing of $\Xi_N$ if 
\begin{enumerate}
    \item for every $t \in [N]$, and $j \in [J]$, $(\Lambda_t^j, \Upsilon_t^j)$ is $\bar{O}(t)$-measurable,
    \item for every $f \in \mathcal{F}$, there exists $j \in [J]$, such that, for every $t \in [J]$, $\Lambda_t^j \leq \xi_t(f) \leq \Upsilon_t^j$,
    \item for every $t \in [N]$, $j \in [J]$, $| \Lambda_t^j - \Upsilon_t^j | \leq b$ a.s.,
    \item for every $j \in [J]$,
    \begin{align}
        \frac{1}{N} \sum_{t=1}^N E \left[ ( \Upsilon_t^j - \Lambda_t^j)^2 | \bar{O}(t-1) \right] \leq \epsilon^2.
    \end{align}
\end{enumerate}
We denote $\mathcal{N}_{[\,]}(\epsilon, b, \Xi_N, \bar{O}(N))$ the minimal cardinality of an $(\epsilon, b, \Xi_N, \bar{O}(N))$-bracketing.
\end{definition}

Applied to our problem, observe that the process $\{M_{2,N}(\bar{Q}, \bar{Q}_1) : \bar{Q} \in \mathcal{Q}\}$ is derived from the process 
\begin{align}
    \Xi_N := \left\lbrace \left((D^*(\bar{Q}) -D^*(\bar{Q}_1))(C_o(t), O(t))\right)_{t=1}^N : \bar{Q} \in \bar{\mathcal{Q}} \right\rbrace.
\end{align}
\noindent
Natural questions that arise are (1) how to connect the sequential bracketing entropy of the process $\Xi_N$ to a traditional bracketing entropy measure for the outcome model $\bar{\mathcal{Q}}$, and (2) how to obtain consistency of an estimator $\bar{Q}^*_N$ fitted from sequentially collected data. Answers to both of these questions entail bracketing entropy preservation results that we present in the upcoming subsection, \ref{subsection:bktg_preservation}.

We emphasize that the notion of \textit{sequential covering numbers}, and the corresponding \textit{sequential covering entropy} introduced by \cite{rakhlin2014}, represent a measure of complexity under which one can control martingale processes and obtain equicontinuity results. One motivation for the development of the notion of sequential covering numbers is that results that hold for i.i.d. empirical processes under traditional covering entropy conditions do not hold for martingale processes. Interestingly, while classical covering number conditions cannot be used to control martingale processes, classical bracketing number bounds can usually be turned into sequential bracketing number bounds. Our choice to state results in terms of one measure of sequential complexity rather than the other (or both) is motivated by concision purposes, and also by the fact that we know how to bound bracketing entropy of a certain class of statistical models we find realistic in many applications, as we describe in later subsections.

\subsubsection{Bracketing preservation results}\label{subsection:bktg_preservation}

We formalize the connection between the sequential bracketing entropy of the process $\Xi_N$ to a traditional bracketing entropy measure for the outcome model $\bar{\mathcal{Q}}$ in lemma \ref{lemma:seq_bktg_can_gdt_process_from_bktg_outcome_model} below. In particular, lemma \ref{lemma:seq_bktg_can_gdt_process_from_bktg_outcome_model} bounds the sequential bracketing entropy of the canonical gradient process $\Xi_N$ in terms of the bracketing entropy of the outcome model $\bar{\mathcal{Q}}$ w.r.t. a norm defined below.

\begin{lemma}\label{lemma:seq_bktg_can_gdt_process_from_bktg_outcome_model}
Suppose that assumption \ref{assumption:strong_positivity} holds. Then 
\begin{align}
    \mathcal{N}_{[\,]}(\epsilon, \Xi_N, \bar{O}(N)) \lesssim N_{[\,]}(\epsilon, \bar{\mathcal{Q}}, L_2(P_{g^*,h_N})),
\end{align}
where $P_{g*,h_N}(a,c) = g^*(a\mid c) h_N(c)$, with $h_N$ being the empirical measure $h_N:= N^{-1} \sum_{t=1}^N \delta_{C_o(t)}$.
\end{lemma}

\begin{proof}
Suppose $\mathcal{B} = \{ (\lambda_j, \upsilon_j) : j \in [J] \}$ is an $\epsilon$-bracketing in $L_2(P_{g^*, h_N})$ norm of $\bar{\mathcal{Q}}$. Let $\bar{Q} \in \mathcal{Q}$. There exists $j \in [J]$ such that $\lambda_j \leq \bar{Q} \leq \upsilon_j$. Without loss of generality, we can suppose that $0 \leq \lambda_j \leq \upsilon_j \leq 1$, since the bracket $(\lambda_j \vee 0, \upsilon_j \wedge 1)$ brackets the same functions of $\bar{\mathcal{Q}}$ as $(\lambda_j, \upsilon_j)$, as every element of $\bar{\mathcal{Q}}$ has range in $[0,1]$. We have that
\begin{align}
    D^*(\bar{Q}) - D^*(\bar{Q}_1) = \frac{g^*}{g_{0,t}} (\bar{Q}_1 - \bar{Q}) + \sum_{a=1}^2 g^*(a \mid \cdot) (\bar{Q} - \bar{Q}_1))(a, \cdot).
\end{align}
Denoting 
\begin{align}
    \Lambda_t^j :=& \frac{g^*}{g_{0,t}} (\bar{Q}_1 - \upsilon_j) + \sum_{a=1}^2 g^*(a \mid \cdot) (\lambda_j - \bar{Q}_1))(a, \cdot), \\
    \text{and }\Upsilon_t^j :=& \frac{g^*}{g_{0,t}} (\bar{Q}_1 - \lambda_j) + \sum_{a=1}^2 g^*(a \mid \cdot) (\upsilon_j - \bar{Q}_1))(a, \cdot),
\end{align}
we have that 
\begin{align}
    \Lambda_t^j \leq (D^*(\bar{Q}, g_{0,t}) - D^*(\bar{Q}, g_{0,t})(C_o(t), O(t)) \leq \Upsilon_t^j.
\end{align}
We now check the size of the sequential bracket $(\Lambda_t^j, \Upsilon_t^j)_{t=1}^N$. We have that
\begin{align}
    &\frac{1}{N} \sum_{t=1}^N E_{Q_0, g}\left[ (\Upsilon_t^j - \Lambda_t^j)^2 \mid \bar{O}(t-1) \right] \\
    =& \ \frac{1}{N} \sum_{t=1}^N E_{Q_0, g_0} \left[ \left\lbrace \frac{g^*}{g_{0,t}}(\upsilon_j - \lambda_j)(A(t), C_o(t)) + \sum_{a=1}^2 (g^* (\upsilon_j - \lambda_j))(a, C_o(t))   \right\rbrace^2 \mid C_o(t) \right] \\
    \leq& \ \frac{2}{N} \sum_{t=1}^N E_{Q_0,g_0} \left[ \left( \frac{g^*}{g_{0,t}} \right)^2 (\upsilon_j - \lambda_j)^2(A(t), C_o(t)) \mid C_o(t) \right] \\
    + & \ E_{Q_0, g^*}\left[ (\upsilon_j - \lambda_j)(A(t), C_o(t))\right]^2 \\
    \leq& \ \frac{4 \delta^{-1}}{N} \sum_{t=1}^N E_{Q_0, g^*} \left[ (\upsilon_j - \lambda_j)^2(A(t), C_o(t)) \mid C_o(t) \right] \\
    =& \ 4 \delta^{-1} \|\upsilon_j - \lambda_j\|_{2, P_{g^*, h_N}}^2 \\
    \leq& \ 4 \delta^{-1} \epsilon^2,
\end{align}
where we have used assumption \ref{assumption:strong_positivity} and Jensen's inequality in the fourth line above. From assumption \ref{assumption:strong_positivity}, it is also immediate to check that $| \Upsilon_t^j - \Lambda_t^j| \leq 2 \delta^{-1}$.

So far, we have proven that one can construct a $(2\delta^{-1/2} \epsilon, 2 \delta^{-1}, \bar{O}(N))$ bracketing of $\Xi_N$ from an $\epsilon$-bracketing in $L_2( P_{g^*, h_N})$ norm of $\bar{\mathcal{Q}}$. Treating $\delta$ as a constant, this implies that $\log N_{[\,]}(\epsilon, 2 \delta^{-1}, \Xi_N, \bar{O}(N)) \lesssim \log N_{[\,]}(\epsilon, \bar{\mathcal{Q}}, L_2(P_{g^*, h_N}))$.
\end{proof}

When proving consistency and convergence rate results for the outcome model estimator $\bar{Q}^*_N$, we need bounds on the sequential bracketing entropy of the following martingale process:
\begin{align}
    \mathcal{L}_N := \left\lbrace \left(\ell_t(\bar{Q})(C_o(t), O(t)) \right)_{t=1}^N : \bar{Q} \in \bar{\mathcal{Q}} \right\rbrace,
\end{align}
where $\ell_t(\bar{Q})(c, o) := (g^*(a \mid c) / g_{0,t}(a \mid c)) (\ell(\bar{Q})(o) - \ell(\bar{Q}_1)(o))$, with $\ell(f)$ denoting a loss function. We refer to $\mathcal{L}_N$ as an \textit{inverse propensity weighted loss process}. Lemma 4 in \cite{bibaut2019fast} provides conditions, which hold for most common loss functions, under which the bracketing entropy of the loss class $\{  \ell(f)(\bar{Q}) : \bar{Q} \in \bar{\mathcal{Q}}\}$ is dominated up to a constant by the bracketing entropy of $\bar{\mathcal{Q}}$. As a direct corollary of this lemma, we state the following result on the sequential bracketing entropy of the process $\mathcal{L}_N$; we refer to \cite{bibaut2019fast} for examples of common settings where assumption \ref{assumption:loss} is satisfied.

\begin{assumption}\label{assumption:loss}
The loss function can be written as $\ell(\bar{Q})(c,a,y) = \widetilde{\ell}(\bar{Q}(c,a),y)$, where $\widetilde{\ell}$ satisfies the following conditions:
\begin{itemize}
    \item for all $f$, $c$, $a$, $y \mapsto \widetilde{\ell}(\bar{Q}(c, a), y)$ is unimodal,
    \item for any $y$, $u \mapsto \widetilde{\ell}(u, y)$ is $L$-Lispchitz, for $L = O(1)$.
\end{itemize}
\end{assumption}

\begin{lemma}[Sequential bracketing entropy of loss process]\label{lemma:seq_bktg_loss_process}
Suppose that assumptions \ref{assumption:loss} and  \ref{assumption:strong_positivity} hold. Then
\begin{align}
    \mathcal{N}_{[\,]}(\epsilon, \mathcal{L}_N, \bar{O}(N)) \lesssim N_{[\,]}(\epsilon, \mathcal{\bar{Q}}, L_2(P_{g^*, h_N})).
\end{align}
\end{lemma}

\subsubsection{Convergence rate of sequentially fitted outcome model estimators}

In this subsection, we give convergence guarantees for outcome model estimators $\bar{Q}_N$, and their targeted counterpart $\bar{Q}_N^*$, fitted on sequentially collected data. We first give convergence rate guarantees for empirical risk minizers $\bar{Q}_N$ over a class $\bar{\mathcal{Q}}$, in terms of the bracketing entropy in $L_2(P_{g^*, h_N})$-norm of $\bar{\mathcal{Q}}$. 

As briefly defined in section \ref{tmle}, let $\ell=L$ be a loss function for the outcome regression such that, for every $\bar{Q}: \mathcal{C} \times \mathcal{A} \to [0,1]$, we have that
\begin{align}
    \bar{Q}_0 \in \arg\min_{\bar{Q} \text{-measurable}} P_{Q_0, g^*, h_N} \ell(\bar{Q}).
\end{align}
We denote $R_{0,N}(\bar{Q}) := P_{Q_0, g^*, h_N} \ell(\bar{Q})$ as the population risk; we note that this population risk is equal to the average across $t$ of the conditional risks $P_{Q_0 g^*, h_N} \ell(\bar{Q})$ given $C_o(t)$. Let $\bar{Q}^*$ be a minimizer of $R_{0,N}(\bar{Q})$ over $\bar{\mathcal{Q}}$. We further define the empirical risk as
\begin{align}
    \widehat{R}_N(\bar{Q}) := \frac{1}{N} \sum_{t=1}^N \frac{g^*}{g_{0,t}}(A(t) \mid C_o(t)) \ell(\bar{Q})(C_o(t), O(t)).
\end{align}
Note that the empirical risk minimizer over $\bar{\mathcal{Q}}$ is any minimizer over $\bar{\mathcal{Q}}$ of $\widehat{R}_N(\bar{Q})$; as such, we use importance sampling weighting factor $g^* / g_{0,t}$ in front of each term $\ell(\bar{Q})(C_o(t), O(t))$. This choice is motivated by the fact that we want convergence rates guarantees for $\bar{Q}_N$ in $L_2(P_{g^*, h_N})$, as is natural to control the size of the sequential brackets of the canonical gradient process $\Xi_N$ in terms of the size of brackets of $\bar{\mathcal{Q}}$ in $L_2(P_{g^*, h_N})$ norm (see lemma \ref{lemma:seq_bktg_can_gdt_process_from_bktg_outcome_model}). In the following, we state the entropy condition and additional assumptions on the loss function. 

\begin{assumption}[Entropy of the outcome model]\label{assumption:entropy_outcome_model}
Suppose that there exists $p > 0$ such that
\begin{align}
    \log (1+N_{[\,]}(\epsilon, \bar{\mathcal{Q}}, L_2(P_{g^*, h_N}))) \leq \epsilon^{-p}.
\end{align}
\end{assumption}
\noindent

\begin{assumption}[Variance bound for the loss]\label{assumption:loss_var_bound}
Suppose that 
\begin{align}
    \|\ell(\bar{Q}) - \ell(\bar{Q}^*) \|^2_{2, \bar{Q}_0, g^*, h_N} \lesssim R_{0,N}(\bar{Q}) - R_{0,N}(\bar{Q}^*)
\end{align}
for all $\bar{Q} \in \bar{\mathcal{Q}}$.
\end{assumption}

\begin{assumption}[Excess risk dominates $L_2$ norm]\label{assumption:excess_risk_dominates_l2_norm}
Suppose that 
\begin{align*}
    \|\bar{Q} - \bar{Q}^*\|_{2, g^*, h_N} \lesssim R_{0,N}(\bar{Q}) - R_{0,N}(\bar{Q}^*).
\end{align*}
\end{assumption}

\begin{theorem}\label{thm:convergence_rate_seq_ERMs}
Consider an empirical risk minimizer $\bar{Q}_N$ over $\bar{\mathcal{Q}}$, and a population minimizer $\bar{Q}^*$, as defined above. Suppose that assumptions \ref{assumption:entropy_outcome_model}, \ref{assumption:loss_var_bound}, \ref{assumption:excess_risk_dominates_l2_norm}, and assumption \ref{assumption:loss} hold. Then,
\begin{align}
    \|\bar{Q}_N - \bar{Q}^* \|_{2, g^*, h_N} = 
    \begin{cases}
    O_P(N^{-\frac{1}{1+p/2}}) & \text{ if } p < 2, \\
    O_P(N^{-\frac{1}{p}}) & \text{ if } p > 2.
    \end{cases}
\end{align}
\end{theorem}

\begin{proof}[Proof]
Consider the process $\mathcal{L}_N$ defined in subsection \ref{subsection:bktg_preservation}. We define $M_{0,N}(\bar{Q}, \bar{Q}^*)$ and $\widehat{M}_N(\bar{Q}, \bar{Q}^*)$ as population and empirical risk differences
\begin{align}
    M_{0,N}(\bar{Q}, \bar{Q}^*) := R_{0,N}(\bar{Q}) - R_{0,N}(\bar{Q^*}) \qquad \text{and} \qquad \widehat{M}_N(\bar{Q}, \bar{Q}^*) := \widehat{R}_N(\bar{Q}) - \widehat{R}_N(\bar{Q}^*). 
\end{align}
Let 
\begin{align}
    \sigma_N^2(\bar{Q}, \bar{Q}^*) := \frac{1}{N} \sum_{t=1}^N E \left[ \left(\frac{g^*}{g_{0,t}}(A(t) \mid C_o(t)) ( \ell(\bar{Q}) - \ell(\bar{Q}^*) )(C_o(t), O(t)) \right)^2\mid C_o(t) \right].
\end{align}
The quantity $\sigma_N(\bar{Q}, \bar{Q}^*)$ can be seen as a sequential equivalent of an $L_2$ norm for the process $\{(g^*/g_{0,t})(A(t) \mid C_o(t)) (\ell(\bar{Q}) - \ell(\bar{Q}^*))(C_o(t), O(t))\}_{t=1}^N $. From assumption \ref{assumption:strong_positivity}, we have that $\sup_{t \geq 1} \|(g^*/g_{0,t}) (\ell(\bar{Q}) - \ell(\bar{Q}^*)) \|_\infty = O(1)$. From theorem A.4 in \cite{vanHandel2010}, with probability at least $1- 2 e^{-x}$, we have that
\begin{align}
    &\sup \left\lbrace M_{0,N}(\bar{Q}, \bar{Q}^*) - \widehat{M}_N(\bar{Q}, \bar{Q}^*) : \bar{Q} \in \bar{\mathcal{Q}}, \ \sigma_N(\bar{Q}) \leq r \right\rbrace\\
    \lesssim & \ r^- + \frac{1}{\sqrt{N}} \int_{r^-}^r \sqrt{\log ( 1 + N_{[\,]}(\epsilon, 1, \mathcal{L}_N, \bar{O}(N)))} d \epsilon  \\
    &+ \ \frac{1}{N}  \log ( 1 + N_{[\,]}(r, 1, \mathcal{L}_N, \bar{O}(N))) + r \sqrt{\frac{x}{N}} + \frac{x}{N}.
\end{align}
From assumption \ref{assumption:strong_positivity}, we have that
\begin{align}
    \sigma_N(\bar{Q}) \lesssim \|\ell(\bar{Q}) - \ell(\bar{Q}^*) \|_{2, g^*,h_N} \lesssim M_{0,N}(\bar{Q}, \bar{Q}^*).
\end{align}
Combined with lemma \ref{lemma:seq_bktg_loss_process}, we have that
\begin{align}
    &\sup \left\lbrace M_{0,N}(\bar{Q}, \bar{Q}^*) - \widehat{M}_N(\bar{Q}, \bar{Q}^*) : \bar{Q} \in \bar{\mathcal{Q}}, \ M_{0,N}(\bar{Q}, \bar{Q}^*) \leq r \right\rbrace \\
    \lesssim & \ r^- + \frac{1}{\sqrt{N}} \int_{r^-}^r \sqrt{\log (1 + N_{[\,]}(\epsilon, \bar{\mathcal{Q}}, L_2(P_{g^*,h_N}))} d\epsilon \\
    &+ \  \frac{1}{N} \log (1 + N_{[\,]}(r, \bar{\mathcal{Q}}, L_2(P_{g^*,h_N})) + r \sqrt{\frac{x}{N}} + \frac{x}{N} \label{eq:pf_rate_ERM_vH_bound}
\end{align}
with probability at least $1 - 2e^{-x}$. In the following, we treat the cases $p < 2$ and $p > 2$ separately.

\paragraph{\textbf{Case p $>$ 2.}} Observe that
\begin{align}
     \| \bar{Q}_N - \bar{Q}^* \|_{2,g^*,h_N} \lesssim & M_{0,N}(\bar{Q}_N, \bar{Q}^*) \\
     =& M_{0,N}(\bar{Q}_N, \bar{Q}^*) - \widehat{M}_N(\bar{Q}_N, \bar{Q}^*) + \widehat{M}_N(\bar{Q}_N, \bar{Q}^*)\\
     \leq & M_{0,N}(\bar{Q}_N, \bar{Q}^*) - \widehat{M}_N(\bar{Q}_N, \bar{Q}^*) \\
     \leq & \sup \left\lbrace M_{0,N}(\bar{Q}, \bar{Q}^*) - \widehat{M}_N(\bar{Q}, \bar{Q}^*) : \bar{Q} \in \bar{\mathcal{Q}}, \ M_{0,N}(\bar{Q}, \bar{Q}^*) \leq r_0 \right\rbrace
\end{align}
where $r_0 := \sup_{\bar{Q} \in \bar{\mathcal{Q}}} M_{0,N}(\bar{Q}, \bar{Q}^*)$. The third line follows from the fact that $Q_N$ minimizes $\widehat{R}_N(\bar{Q})$ over $\bar{\mathcal{Q}}$, wich implies that $\widehat{M}_N(\bar{Q}_N, \bar{Q}^*) \leq 0$. We now use equation \eqref{eq:pf_rate_ERM_vH_bound} to bound the last line of the inequality. From assumption \ref{assumption:strong_positivity}, we know that $r_0 = O(1)$. Using the entropy bound from assumption  \ref{assumption:entropy_outcome_model} and minimizing the right hand side of \eqref{eq:pf_rate_ERM_vH_bound} w.r.t. $r^-$, we obtain that, with probability at least $1 - 2e^{-x}$,
\begin{align}
    \| \bar{Q}_N - \bar{Q}^* \|^2_{2,g^*,h_N} \lesssim N^{-2/p} + \frac{x}{\sqrt{N}} + \frac{x}{N},
\end{align}
which, by picking $x$ appropriately, then implies that $\|\bar{Q}_N - \bar{Q}^*\|_{2,g^*,h_N} = O_P(N^{-1/p})$.

\paragraph{\textbf{Case p $<$ 2.}} Starting from the bound \eqref{eq:pf_rate_ERM_vH_bound}, via some algebra and by taking an integral, we obtain
\begin{align}
    &E_{P_0} \left[ \sup \left\lbrace M_{0,N}(\bar{Q}, \bar{Q}^*) - \widehat{M}_N(\bar{Q}, \bar{Q}^*) : \bar{Q} \in \bar{\mathcal{Q}}, \ M_{0,N}(\bar{Q}, \bar{Q}^*) \leq r \right\rbrace \right] \\
    \lesssim & \ r^- + \frac{1}{\sqrt{N}} \left( r + \int_{r^-}^r \sqrt{\log (1 + N_{[\,]}(\epsilon, \bar{\mathcal{Q}}, L_2(P_{g^*,h_N}))} d\epsilon \right) \\
    &+ \ \frac{1}{N} \left(r + \log (1 + N_{[\,]}(r, \bar{\mathcal{Q}}, L_2(P_{g^*,h_N})) \right).
\end{align}
Let $r^- = 0$. By using the entropy bound from assumption \ref{assumption:entropy_outcome_model}, we obtain that
\begin{align}
    &E_{P_0} \left[ \sup \left\lbrace M_{0,N}(\bar{Q}, \bar{Q}^*) - \widehat{M}_N(\bar{Q}, \bar{Q}^*) : \bar{Q} \in \bar{\mathcal{Q}}, \ M_{0,N}(\bar{Q}, \bar{Q}^*) \leq r \right\rbrace \right] \\
    &\lesssim \frac{1}{\sqrt{N}} r^{1 - p/2} \left(1 + \frac{r^{1-p/2}}{r^2 \sqrt{N}} \right).
\end{align}
Theorem 3.4.1 in \cite{vanderVaart&Wellner96} then implies that
\begin{align}
    M_{0,N}(\bar{Q}_N, \bar{Q}^*) = O_P(N^{-\frac{2}{1 + p/2}}), \text{ and therefore } \| \bar{Q}_N - \bar{Q}^* \|_{2,g^*,h_N} = O_P(N^{-\frac{1}{1 + p/2}}).
\end{align}
\end{proof}

\subsection{Outcome model classes}\label{halmore}

Now that we know how to characterize the sequential bracketing entropy of $\Xi_N$ and $\mathcal{L}_N$ in terms of the bracketing entropy w.r.t. the norm $L_2(P_{Q_0, h_{C,N}})$ of the outcome model $\mathcal{Q}$, we look at specific function classes $\mathcal{Q}$ for which we know how to bound the latter.

\subsubsection{Holder class}
Consider functions over a certain domain $\mathcal{X}$; in our setting we note that $\mathcal{X} = \mathcal{C} \times \mathcal{O}$. Suppose that $\mathrm{dim}(\mathcal{X}) = d$. 
We denote $H(\beta, M)$ the class of functions over a certain domain $\mathcal{X}$ , such that, for any $x, y \in \mathcal{X}$, and any non-negative integers $\beta_1,\ldots, \beta_d $ such that $\beta_1 + \ldots + \beta_d = \left\lfloor \beta \right\rfloor$, 
\begin{align}
    \left\lvert \frac{\partial^{\left\lfloor \beta \right\rfloor} f } { \partial x_1^{\beta_1} \ldots \partial x_d^{\beta_d}}(x) - \frac{\partial^{\left\lfloor \beta \right\rfloor} f } { \partial x_1^{\beta_1} \ldots \partial x_d^{\beta_d}}(y) \right\rvert\leq  M \|x - y\|.
\end{align}
The bracketing entropy w.r.t. the uniform norm $\|\cdot\|$ of such a class satisfies
\begin{align}
    \log N_{[\,]}(\epsilon, H(\beta, M), \|\cdot\|) \lesssim \epsilon^{-d/\beta}.
\end{align}
For more detail, we refer the interested reader to, for example, chapter 2.7 in \cite{vanderVaart&Wellner96}. As such, our Donsker condition \ref{assumption:Donsker_condition_EIF_process} is satisfied for $\beta > d / 2$. Nevertheless, we caution that assuming that the outcome model lies in a Holder class of differentiability order $\beta > d/2$ might be an overly restrictive assumption.

\subsubsection{HAL class}
A class of functions that is much richer that the previous Holder classes is the class of cadlag functions with bounded sectional variation norm - also referred to as Hardy-Krause variation. We refer to this class as the Highly Adaptive Lasso class (HAL class), as it is the class in which the estimator, introduced in \cite{vanderlaan2017hal}, takes values. The Highly Adaptive Lasso class is particularly attractive in i.i.d. settings for various reasons, which we enumerate next. (1) Unlike Holder classes, it doesn't make local smoothness assumptions. Rather it only restricts a global measure of irregularity, the sectional variation norm, thereby allowing for functions to be differentially smooth/variable depending on the area of the domain. (2) Emprical risk minimzers over the HAL class were shown to be competitive with the best supervised machine learning algorithms, including Gradient Boosting Machines and Random Forests. (3) We know how to bound both the uniform metric entropy and the bracketing entropy of these classes of functions. These bounds show that the corresponding entropy integrals are bounded, which imply that the HAL class is Donsker. In particular, \cite{bibaut2019fast} provide a bound on the bracketing entropy w.r.t. $L_r(P)$, for $r \in [1, \infty)$, for probability distribution that have bounded Radon-Nikodym derivative w.r.t. the Lebesgue measure, that is $dP / d\mu \leq C$. \cite{bibaut2019fast} use this bracketing entropy bound to prove the rate of convergence $O(N^{-1/3} (\log N)^{2d-1})$.

Unfortunately, to bound the sequential bracketing entropies of $\Xi_N$ and of $\mathcal{L}_N$ we would need a bracketing entropy bound w.r.t. $L_2(P_{Q_0, h_{C,N}})$, which, owing to the fact that $h_{C,N}$ is a discrete measure, does not have bounded Radon-Nikodym derivative w.r.t. the Lebesgue measure over $\mathcal{C} \times \mathcal{O}$. Under the assumption \ref{assumption:context_marginal_convergence} on the convergence of the marginals of $(C_o(t))$ to a limit law (we shall denote it $h_{\infty}$), we have that $h_{C,N} \xrightarrow{d} h_{\infty}$, which can reasonably be a continuous measure dominated by the Lebesgue measure.  By convergence in distribution of $h_{C,N}t$ to $h_{\infty}$, we have at least that the size of brackets w.r.t. $h_{C,N}$ converges to the size of brackets under $h_{\infty}$. If this convergence were uniform over bracketings of $\bar{\mathcal{Q}}$, and that $dh_{\infty} / d\mu \leq C$, then we would have that $N_{[\,]}(\epsilon, \bar{\mathcal{Q}}, L_2(P_{Q_0, h_{C,N}})) \lesssim N_{[\,]}(\epsilon, \bar{\mathcal{Q}}, L_2(\mu)$. Proving the uniformity over bracket seems to be a relatively tough theoretical endeavor, and we leave it to future research.

\subsubsection{A modified HAL class}
Given the difficulty in bounding $N_{[\,]}(\epsilon, \bar{\mathcal{Q}}, L_2(P_{Q_0, g^*, h_N}))$ for the HAL, class, we consider a modified HAL class in the case where $\mathcal{C}$ is discrete, that is $\mathcal{C} = \{c_1,\ldots, c_J \}$. We define the modified class as the set of functions $f : \mathcal{C} \times \mathcal{O} \to \mathbb{R}$ such that, for every $c \in \mathcal{C}$, $o \mapsto f(c,o)$ is cadlag with sectional variation norm smaller than $M_1$. It is straightforward to show that the bracketing entropy of such a class $\mathcal{F}$ is bounded as follows:
\begin{align}
    \log N_{[\,]}(\epsilon, \mathcal{F}, L_2(P_{Q_0,g^*,h_N})) \lesssim |\mathcal{C}| \epsilon^{-1} (\log (1/\epsilon))^{2(\dim(\mathcal{O}) - 1)}.
\end{align}

\subsection{Proof of Lemma \ref{lemma:cond_var_stab_under_mixing}}\label{lemma:cond_var_stab_under_mixing_proof}

The proof of lemma \ref{lemma:cond_var_stab_under_mixing} relies on lemma \ref{lemma:cond_var_Lipschitz_wrt_design}, which we present and prove in the following.

\begin{lemma}\label{lemma:cond_var_Lipschitz_wrt_design}
Denote, for any fixed $g$,
\begin{align}
    f(g)(c) := \Var_{Q_0,g_0} \left(D^*(Q_1, g)(C_o(t), O(t)) \mid C_o(t) = c \right).
\end{align}
Suppose that assumption \ref{assumption:strong_positivity} holds, and
let $g$ be fixed given $C_o(t)$. Suppose that the strong positivity assumption holds for $g$, that is $g \geq \delta$, for the same $\delta$ an in assumption \ref{assumption:strong_positivity}. Then,
\begin{align}
    \left\lvert E_{h_{0,t}} [  f(g)(C_o(t)) ]  -  E_{h_{0,t}} [f(g_{0,t})(C_o(t))]\right\rvert \leq 4 \delta^{-3} \|g - g_{0,t} \|_{1, P_{g^*, h_{0,t}}}.
\end{align}
\end{lemma}

\begin{proof}
Observe that $D^*(Q_1,g)(c,o)$ can be decomposed as 
\begin{align*}
    D^*(Q_1, g)(c,o) = D^*_1(\bar{Q}_1,g) + D^*_2(Q_1, c),
\end{align*}
with 
\begin{align}
    D^*_1(\bar{Q}_1, g)(c,o) :=& \frac{g^*(a\mid c)}{g(a\mid c)} \left( y - \bar{Q}_1(a,c) \right), \\
\text{and }   D_2^*(Q_1)(c) :=& \sum_{a=1}^2 g^*(a' \mid c) \bar{Q}_1(a',c) - \Psi(Q_1).
\end{align}
As $D^*_2(Q_1)(C_o(t))$ is constant given $C_o(t)$, we have that
\begin{align}
    f(g)(c) = \Var_{Q_0, g_0} \left(D^*_1(\bar{Q}_1, g)(C_o(t), O(t)) \mid C_o(t) \right).
\end{align}
\noindent
In the following, the dependence of the canonical gradient on $(C_o(t),O(t))$ is implied, but suppressed in the notation. For any $g$, 
\begin{align}
      &\left\lvert E_{h_{0,t}} [  f(g)(C_o(t)) ]  -  E_{h_{0,t}} [f(g_{0,t})(C_o(t))]\right\rvert  \\
      =& \left\lvert E_{h_{0,t}} \left[ E_{Q_0, g_0} \left[ (D_1^*(\bar{Q}_1, g))^2 \mid C_o(t)) \right] \right]
      - \left[ E_{Q_0, g_0} \left[ (D_1^*(\bar{Q}_1, g_{0,t}))^2 \mid C_o(t)) \right] \right] \right. \\
      & - \left. E_{h_{0,t}} \left[ E_{Q_0, g_0} \left[D^*_1(\bar{Q}_1, g) \mid C_o(t) \right]^2 -  E_{Q_0, g_0} \left[D^*_1(\bar{Q}_1, g_{0,t}) \mid C_o(t) \right]^2 \right]\right\rvert \\
      \leq & \ E_{Q_0, g_0} \left[ \left\lvert (D^*_1(\bar{Q}_1, g))^2 -(D^*_1(\bar{Q}_1, g_{0,t}))^2 \right\rvert \right] \\
      &+ \left(\| D^*_1(\bar{Q}_1, g) \|_\infty + \|D^*_1(\bar{Q}_1, g_{0,t}) \|_\infty \right) \times E_{Q_0, g_0} \left[ \left\lvert D^*_1(\bar{Q}_1, g_1) - D_1^*(\bar{Q}_1, g_{0,t}) \right\rvert (C_o(t), O(t))) \right].
\end{align}
\noindent
We start with the analysis of the first term. In particular, we have that
\begin{align}
    & E_{Q_0, g_0} \left[ \left\lvert (D^*_1(\bar{Q}_1, g))^2 -(D^*_1(\bar{Q}_1, g_{0,t}))^2 (C_o(t), O(t)) \right\rvert \right] \\
    =& \ E_{Q_0, g_0} \left[ \left(Y(t) - \bar{Q}_1(A(t), C_o(t)) \right)^2 \left\lvert \left(\frac{g^*}{g}\right)^2(A(t) \mid C_o(t)) - \left(\frac{g^*}{g_{0,t}}\right)^2 (A(t) \mid C_o(t))\right\rvert\right]  \\
    =& \ E_{Q_0, g_0} \left[ \left( Y(t) - \bar{Q}_1(A(t) \mid C_o(t))\right)^2 \left( \frac{g^*}{g_{0,t}} \right)^2 (A(t) \mid C_o(t))\left\lvert 1 - \left( \frac{g}{g_{0,t}}\right)^2(A(t) \mid C_o(t))\right\rvert \right] \\
    \leq & \ \delta^{-1} E_{Q_0, g^*} \left[ \left( Y(t) - \bar{Q}_1(A(t), C_o(t)) \right)^2 \left\lvert \left(\frac{g_{0,t}^2 - g^2}{g^2}\right)(A(t), C_o(t)) \right\rvert\right] \\
    \leq & \ \delta^{-3} E_{Q_0, g^*} \left[ \left\lvert (g_{0,t}^2 - g^2)(A(t) \mid C_o(t)) \right\rvert\right] \\
    \leq & \ 2 \delta^{-3} E_{Q_0, g^*} \left[ \left\lvert (g_{0,t} - g)(A(t) \mid C_o(t))\right\rvert \right] \\
    =& 2 \delta^{-3} \| g_{0,t} - g \|_{1, P_{g^*, h_{0,t}}}.
\end{align}
\noindent
We now turn to the second term. It follows that
\begin{align}
    E_{Q_0, g_0} \left[ \left\lvert D_1^*(\bar{Q}, g)(C_o(t), O(t)) - D_1^*(C_o(t), O(t))  \right\rvert \right]
    =& \ E_{Q_0, g_0} \left[ \left( \bar{Q}_1 \left\lvert \frac{g^*}{g} - \frac{g^*}{g_{0,t}}\right\rvert \right)(C_o(t), O(t)) \right] \\
    \leq & \ \delta^{-1} \| g_{0,t} - g \|_{1, P_{g^*, h_{0,t}}} \\
    \leq & \ \delta^{-1} E_{Q_0, g^*} \left[ \left\lvert (g_{0,t} - g)(A(t) \mid C_o(t)) \right\rvert\right] \\
    =& \ \delta^{-1} \|g_{0,t} - g \|_{1, P_{g^*, h_{0,t}}}.
\end{align}
As $\|D_1^*(\bar{Q}_1, g)\|_\infty \leq \delta^{-1}$ and $\|D_1^*(\bar{Q}_1, g_{0,t} \|_\infty \delta^{-1}$, we therefore have that 
\begin{align}
    E_{Q_0, g_0} \left[ \left\lvert f(g)(C_o(t)) - f(g_{0,t})(C_o(t)) \right\rvert\right] \leq 4 \delta^{-3} \| g - g_{0,t} \|_{1, P_{g^*, h_{0,t}}}.
\end{align}
\end{proof}

\begin{manuallemma}{2}Suppose that assumptions \ref{assumption:strong_positivity}, \ref{assumption:context_marginal_convergence} and \ref{assumption:rho_mixing} hold. Then assumption \ref{assumption:cond_variances} holds.
\end{manuallemma}
\noindent

\begin{proof}

We use the notation of lemma \ref{lemma:cond_var_Lipschitz_wrt_design} in this proof. We have that
\begin{align}
    &\frac{1}{N} \sum_{t=1}^N f(g_{0,t})(C_o(t)) - E \left[ f(g_\infty) (C_\infty) \right] \\
    =& \ \frac{1}{N} \sum_{t=1}^N f(g_{0,t})(C_o(t)) - E \left[ f(g_{0,t}(C_o(t)) \right] \label{eq:pf_stab_var_AN} \\
    &+ \ \frac{1}{N} \sum_{t=1}^N E \left[ f(g_{0,t})(C_o(t))\right]- E\left[ f(g_\infty)(C_o(t))\right] \label{eq:pf_stab_var_BN}\\
    &+ \ \frac{1}{N} \sum_{t=1}^N E \left[ f(g_\infty)(C_o(t))\right] - E \left[ f(g_\infty)(C_\infty) \right]\label{eq:pf_stab_var_CN}.
\end{align}
Denote $A_N$, $B_N$ and $C_N$ the quantities in lines \eqref{eq:pf_stab_var_AN}, \eqref{eq:pf_stab_var_BN} and \eqref{eq:pf_stab_var_CN} above; we start with bounding $A_N$. In particular, we have that 
\begin{align}
    \Var\left( A_N \right)=& \frac{1}{N^2} \sum_{t=1}^N \Var(f(g_{0,t})(C_o(t))) \\
    &+ \ \frac{2}{N^2} \sum_{t=1}^N \sum_{s=1}^{N-t} \Cov(f(g_{0,t})(C_o(t)), f(g_{0,t+s})(C_o(t+s))) \\
    =& \ \frac{1}{N^2} \sum_{t=1}^N \Var(f(g_{0,t})(C_o(t))) \\
    &+ \ \frac{2}{N^2} \sum_{t=1}^N \sum_{s=1}^{N-t} \sqrt{\Var(f(g_{0,t})(C_o(t))) \Var(f(g_{0,t+s})(C_o(t+s)))} \rho(C_o(t), C_o(t+s)).
\end{align}
From assumption \ref{assumption:strong_positivity}, $\|f(g_{0,t})\|_\infty \leq \delta^{-1}$ for every $t$. Therefore, 
\begin{align}
    \Var\left(  \frac{1}{N} \sum_{t=1}^N f(g_{0,t})(C_o(t)) \right) 
    \leq & \ \frac{\delta^{-2}}{N^2} \left( N  + 2  N \sup_{t \geq 1} \sum_{s=1}^N \rho(C_o(t), C_o(t+s)) \right) \\
    \leq & \ \frac{\delta^{-2}}{N^2} \left( N + o(N^2)) \right) \\
    =& \ o(1),
\end{align}
where we have used assumption \ref{assumption:rho_mixing} in the last line of the above inequality. Therefore, from Chebyshev's inequality, 
\begin{align}
    \frac{1}{N} \sum_{t=1}^N f(g_{0,t})(C_o(t)) - \frac{1}{N} \sum_{t=1}^N E[f(g_{0,t})(C_o(t))] = o_P(1). \label{eq:vanishing_variance_of_cond_variances}
\end{align}
\noindent
We now turn to $B_N$. From lemma \ref{lemma:cond_var_Lipschitz_wrt_design}, we have that
\begin{align}
    B_N \leq 4 \delta^{-3} \frac{1}{N} \sum_{t=1}^N \|g_{0,t} - g_\infty \|_{1, P_{g^*, h_{0,t}}}.
\end{align}
From assumption \ref{assumption:design_stabilization} and Cesaro's lemma for deterministic sequences of real numbers, $B_N = o(1)$. Finally, from assumption \ref{assumption:context_marginal_convergence} and Cesaro's lemma, $C_N = o(1)$. Denoting $\sigma_0^2 := E[f(C_\infty)]$, the above inequality and \eqref{eq:vanishing_variance_of_cond_variances} yield the wished claim.
\end{proof}

\subsection{Proof of Theorem \ref{thm:TMLE_decompositon}}

\begin{manualtheorem}{2}
For any $\bar{Q}_1 \in \bar{\mathcal{Q}}$, the difference between the TMLE and its target decomposes as 
\begin{align}
    \bar{\Psi}(\bar{Q}_N^*) - \bar{\Psi}(\bar{Q}_0) = M_{1,N}(\bar{Q}_1) + M_{2,N}(\bar{Q}_N^*, \bar{Q}_1),
\end{align}
with 
\begin{align}
    M_{1,N}(\bar{Q}_1) =& \frac{1}{N} \sum_{t=1}^N D^*(\bar{Q}_1) (C_o(t), O(t)) - P_{0, C_o(t)} D^*(\bar{Q}_1), \\
    M_{2,N}(\bar{Q}^*_N, \bar{Q}_1) =& \frac{1}{N} \sum_{t=1}^N (\delta_{C_o(t), O(t)} - P_{0, C_o(t)}) (D^*(\bar{Q}^*_N) - D^*(\bar{Q}_1)).
\end{align}
\end{manualtheorem}

\begin{proof}[Proof] We recall the first order expansion of $\Psi_{C_o(t)}$ given by Theorem \ref{thm:can_gdt_and_1st_order_exp}, and defined as
\begin{align}
    \Psi_{C_o(t)}(\bar{Q}) - \Psi_{C_o(t)}(\bar{Q}_0) = - P_{0,C_o(t)} D^*(\bar{Q}, g)(C_o(t), O(t)) + R(\bar{Q},\bar{Q}_0, g, g_{0,t}).
\end{align}
We note that, since in an adaptive trial the treatment mechanism is controlled, we have that $g = g_0$. Additionally, by definition, TMLE procedure yields $\bar{Q}^*_N$ such that 
\begin{align}
    \frac{1}{N} \sum_{t=1}^N D^*(\bar{Q}^*_N)(C_o(t), O(t)) = 0.
\end{align}
Combined, it follows that
\begin{align}
    \bar{\Psi}(\bar{Q}^*_N) - \bar{\Psi}(\bar{Q}_0) =& \frac{1}{N} \sum_{t=1}^N D^*(\bar{Q}^*_N)(C_o(t), O(t)) - P_{0,C_o(t)} D^*(\bar{Q}^*_N)(C_o(t), O(t)) 
\end{align}
which, by adding and subtracting  $N^{-1}\sum_{t=1}^N (\delta_{C_o(t), O(t)} - P_{0,C_o(t)}) D^*(\bar{Q}_1)$, implies the wished decomposition.
\end{proof}

\subsection{Proof of Theorem \ref{thm:equicontinuity_M2N}}

\begin{manualtheorem}{4}
Consider the process $\Xi_N$ defined in equation \eqref{eq:def_can_gdt_process}. Suppose that assumptions \ref{assumption:strong_positivity}, \ref{assumption:Donsker_condition_EIF_process} and \ref{assumption:convergence_outcome_model} hold. Then $M_{2,N}(\bar{Q}^*_N, \bar{Q}_1) = o_P(N^{-1/2})$.
\end{manualtheorem}

\begin{proof}
We want to show that for any $\epsilon, \delta > 0$, there exist $N_0$ such that for any $N \geq N_0$
\begin{align}
    P_0 \left[ \sqrt{N} M_{2,N}(\bar{Q}^*_N, \bar{Q}_1) \geq \epsilon \right] \leq \delta.
\end{align}
\noindent
Let $\epsilon > 0$ and $\delta > 0$. Define, for any $\bar{Q}$,
\begin{align}
    \sigma^2_N(\bar{Q}, \bar{Q}_1) := \frac{1}{N} \sum_{t=1}^N E_{P_0} \left[ \left(D^*(\bar{Q}, g_{0,t}) - D^*(\bar{Q}_1, g_{0,t}) \right)(C_o(t), O(t))^2  \mid C_o(t) \right].
\end{align}
Under assumption \ref{assumption:strong_positivity}, $\sup_{t \geq 1} \sup_{\bar{Q} \in \bar{\mathcal{Q}}} \|D^*(\bar{Q}, g_{0,t}) - D^*(\bar{Q}_1, g_{0,t})\|_\infty = O(1)$. Theorem A.4 in \cite{vanHandel2010} yields that, with probability at least $
1 - \delta / 2$,
\begin{align}
    & \sqrt{N} \sup \left\lbrace M_{2,N}(\bar{Q}, \bar{Q}_1) : \bar{Q} \in \bar{\mathcal{Q}},\ \sigma_N(\bar{Q},\bar{Q}_1) \leq r \right\rbrace \\
    \lesssim & \ J_{[\,]}(r, 1, \Xi_N, \bar{O}(N)) + \frac{1}{\sqrt{N}} \log( 1 + N_{[\,]}(r, 1, \Xi_N, \bar{O}(N))) + r \sqrt{\log(2/\delta)} + \frac{\log(2/\delta)}{\sqrt{N}}. \label{eq:pf_equicont_vH_bound}
\end{align}
As $\epsilon \mapsto \sqrt{ \log( 1 + N_{[\,]}(\epsilon, 1, \Xi_N, \bar{O}(N)))}$ is non-increasing, we have that 
\begin{align}
    \log( 1 + N_{[\,]}(r, 1, \Xi_N, \bar{O}(N))) \leq \frac{\left(J_{[\,]}(r, 1, \Xi_N, \bar{O}(N))\right)^2}{r^2},
\end{align}
and therefore we can bound the right-hand side in \eqref{eq:pf_equicont_vH_bound} with
\begin{align}
    J_{[\,]}(r, 1, \Xi_N, \bar{O}(N)) \left( 1 + \frac{J_{[\,]}(r, 1, \Xi_N, \bar{O}(N))}{\sqrt{N} r^2} \right) + r \sqrt{\log(2/\delta)} + \frac{\log(2/\delta)}{\sqrt{N}}.
\end{align}
From assumption \ref{assumption:Donsker_condition_EIF_process}, there exists $r_0 > 0$ small enough that
\begin{align}
    r_0 \sqrt{\log (2 / \delta} \leq \epsilon / 4 \qquad \text{and} \qquad J_{[\,]}(r_0, 1, \Xi_N, \bar{O}(N)) \leq \epsilon / 4.
\end{align}
We choose $N_1$ such that, for every $N \geq N_1$,
\begin{align}
    \frac{J_{[\,]}(r_0, 1, \Xi_N, \bar{O}(N))}{\sqrt{N} r_0^2} \leq 1 \qquad \text{and} \qquad \frac{\log (2 / \delta)}{\sqrt{N}} \leq \epsilon / 4.
\end{align}
We then have that, for any $N \geq N_1$, with probability at least $1 - \delta / 2$,
\begin{align}
    \sqrt{N} \sup \left\lbrace M_{2,N}(\bar{Q}, \bar{Q}_1) : \bar{Q} \in \bar{\mathcal{Q}},\ \sigma_N(\bar{Q},\bar{Q}_1) \leq r_0 \right\rbrace \leq \epsilon / 4. \label{eq:def_event_E1}
\end{align}
Denote $\mathcal{E}_1(N,r_0)$ the event under which \eqref{eq:def_event_E1} holds. We further ntroduce the event
\begin{align}
    \mathcal{E}_2(N, r_0) := \left\lbrace \sigma_N(\bar{Q}_N^*, \bar{Q}_1) \leq r_0 \right\rbrace.
\end{align}
Under assumption \ref{assumption:strong_positivity}, $\sigma_N(\bar{Q}_N^*, \bar{Q}_1) \lesssim \| \bar{Q}_N^* - \bar{Q}_1 \|_{2,g^*,h_N}$, and from assumption \ref{assumption:convergence_outcome_model}, we have that 
\begin{align*}
    \|\bar{Q}_N^* - \bar{Q}_1 \|_{2,g^*,h_N} = o_P(1).
\end{align*}
Therefore, there exists $N_2$ such that for every $N \geq N_2$, $\mathcal{E}_2(N, r_0)$ holds with probability at least $1-\delta / 2$. We further conclude that for any $N \geq N_0 := \max(N_1, N_2)$, $P_0[ \mathcal{E}_1(N,r_0) \cap \mathcal{E}_2(N,r_0)] \geq 1 - \delta $, and under $\mathcal{E}_1(N,r_0) \cap \mathcal{E}_2(N,r_0)$, it holds that
\begin{align}
    \sqrt{N} M_{2,N}(\bar{Q}_N^*, \bar{Q}_1) \leq \epsilon,
\end{align}
which is the wished claim.
\end{proof}

\newpage
\bibliographystyle{agsm}
\bibliography{single_ts_adapt}

\end{document}